\documentclass[a4paper]{siamart171218}
\usepackage{amsmath}
\usepackage{amsfonts}
\usepackage{amscd}
\usepackage{amssymb}
\usepackage{graphicx}
\usepackage{color}
\usepackage{booktabs}
\usepackage{mathtools}
\usepackage{hyperref}
\usepackage{tikz}
\usetikzlibrary{trees,calc,positioning,decorations.markings,arrows,backgrounds}
\usetikzlibrary{shapes, shapes.geometric, fit}
\usepackage{standalone}
\usepackage{xspace}
\usepackage{stmaryrd}
\usepackage{enumerate}
\usepackage{subcaption}

\newcommand{\todo}[1]{\textcolor{red}{XXX: #1}}

\definecolor{blue}{rgb}{0.2980392156862745, 0.4470588235294118, 0.6901960784313725}
\definecolor{green}{rgb}{0.3333333333333333, 0.6588235294117647, 0.40784313725490196}
\definecolor{red}{rgb}{0.7686274509803922, 0.3058823529411765, 0.3215686274509804}

\DeclareMathOperator{\vcurl}{\mathbf{curl}}
\newcommand{\Ha}{\mathrm{Ha}}

\newcommand{\B}{\mathbf{B}}
\newcommand{\E}{\mathbf{E}}
\renewcommand{\u}{\mathbf{u}}
\renewcommand{\v}{\mathbf{v}}
\newcommand{\F}{\mathbf{F}}
\newcommand{\C}{\mathbf{C}}

\newcommand{\n}{\mathbf{n}}

\newcommand{\zerov}{\mathbf{0}}
\newcommand{\f}{\mathbf{f}}
\newcommand{\x}{\mathbf{x}}
\newcommand{\V}{\mathbf{V}}
\newcommand{\Rr}{\mathbf{R}}
\newcommand{\W}{\mathbf{W}}

\newcommand{\MM}{\mathcal{M}}
\newcommand{\ZZ}{\mathcal{Z}}
\newcommand{\UU}{\mathcal{U}}
\newcommand{\NN}{\mathcal{N}}
\newcommand{\II}{\mathcal{I}}
\newcommand{\JJ}{\mathcal{J}}
\newcommand{\CC}{\mathcal{C}}
\newcommand{\DD}{\mathcal{D}}
\renewcommand{\AA}{\mathcal{A}}
\newcommand{\BB}{\mathcal{B}}
\newcommand{\GG}{\mathcal{G}}
\newcommand{\FF}{\mathcal{F}}

\renewcommand{\Re}{\ensuremath{\mathrm{Re}}\xspace}
\newcommand{\Rem}{\ensuremath{\mathrm{Re_m}}\xspace}

\newcommand{\grad}{\ensuremath{\nabla}\xspace}
\let\div\relax
\DeclareMathOperator{\div}{div}
\DeclareMathOperator{\curl}{curl}
\newcommand{\scurl}{\curl}

\newcommand{\Ned}{N\'{e}d\'{e}lec\ }

\newcommand{\R}{\mathbb{R}}
\newcommand{\hcurl}{\mathbf{H}(\curl)}
\newcommand{\hdiv}{\mathbf{H}(\div)}
\newcommand{\hdivz}{\mathbf{H}(\div,0)}
\newcommand{\hzcurl}{\mathbf{H}_0(\curl, \Omega)}
\newcommand{\hzdiv}{\mathbf{H}_0(\div, \Omega)}
\newcommand{\Hozv}{\mathbf{H}^1_0(\Omega)}
\newcommand{\Hoz}{H^1_0(\Omega)}

\newcommand{\del}{\partial}

\newcommand{\eps}[1]{\ensuremath{\varepsilon(#1)}}

\newcommand\scalemath[2]{\scalebox{#1}{\mbox{\ensuremath{\displaystyle #2}}}}

\newcommand*{\ldblbrace}{\{\mskip-5mu\{}
\newcommand*{\rdblbrace}{\}\mskip-5mu\}}

\definecolor{darkblue}{rgb}{0.00,0.00,0.55}
\definecolor{black}{rgb}{0.00,0.00,0.00}
\hypersetup{
linkcolor = darkblue,
anchorcolor = darkblue,
citecolor = darkblue,
filecolor = darkblue,
urlcolor = darkblue,
colorlinks  = true,
}

\title{An augmented Lagrangian preconditioner for the
magnetohydrodynamics equations at high Reynolds and coupling numbers
\thanks{Submitted to the editors \today.
	\funding{The first author was supported by the EPSRC Centre for Doctoral Training in Partial Differential Equations: Analysis and Applications, grant EP/L015811/1. The second author was supported by EPSRC grants EP/V001493/1 and EP/R029423/1.}}}

\author{
Fabian~Laakmann\thanks{Mathematical Institute, University of Oxford, Oxford, UK (\email{fabian.laakmann@maths.ox.ac.uk}).}
\and
Patrick~E.~Farrell\thanks{Mathematical Institute, University of Oxford, Oxford, UK (\email{patrick.farrell@maths.ox.ac.uk}).}
\and
Lawrence~Mitchell\thanks{Department of Computer Science, Durham
  University, Durham, UK; Present address: NVIDIA Corporation, Santa
  Clara, USA (\email{lmitchell@nvidia.com}).}}

\headers{Preconditioners for MHD}{F.~Laakmann, P.~E.~Farrell and L.~Mitchell}

\begin{document}
\maketitle

\begin{abstract}
	The magnetohydrodynamics (MHD) equations are generally known to be difficult to solve numerically, due to their highly nonlinear structure and the strong coupling between the electromagnetic and hydrodynamic variables, especially for high Reynolds and coupling numbers.
	In this work, we present a scalable augmented Lagrangian preconditioner for a finite element discretization of the $\B$-$\E$ formulation of the incompressible viscoresistive MHD equations. For stationary problems, our solver achieves robust performance with respect to the Reynolds and coupling numbers in two dimensions and good results in three dimensions. We extend our method to fully implicit methods for time-dependent problems which we solve robustly in both two and three dimensions. Our approach relies on specialized parameter-robust multigrid methods for the hydrodynamic and electromagnetic blocks. The scheme ensures exactly divergence-free approximations of both the velocity and the magnetic field up to solver tolerances.
	We confirm the robustness of our solver by numerical experiments in which we consider fluid and magnetic Reynolds numbers and coupling numbers up to 10,000 for stationary problems and up to 100,000 for transient problems in two and three dimensions.
\end{abstract}



\section{Introduction} \label{sec:introduction}

In this work, we consider the incompressible viscoresistive magnetohydrodynamics (MHD) equations on a simply-connected polytopal domain $\Omega\subset\R^d$, $d \in \{2,3\}$. In the stationary three-dimensional setting, we investigate the formulation
\begin{subequations}
	\label{eq:MHD3d}
	\begin{align}
		- \frac{2}{\Re} \div \eps \u + \u \cdot \nabla \u + \nabla p + S\, \B \times (\E + \u \times \B) &= \f, \label{eq:MHD3d1}\\
		\div \u&=0, \label{eq:MHD3d2}\\
		\E + \u \times \B - \frac{1}{\Rem} \vcurl \B &= \mathbf{0}, \label{eq:MHD3d3}\\
		\vcurl \E &= \mathbf{0}, \label{eq:MHD3d4}\\
		\div \B &= 0. \label{eq:MHD3d5}
	\end{align}
\end{subequations}
Here, $\u:\Omega\to\R^3$ denotes the velocity, $p:\Omega\to\R$ the fluid pressure, $\B:\Omega\to\R^3$ the magnetic field, $\E:\Omega\to\R^3$ the electric field, $\Re$ the fluid Reynolds number, $\Rem$ the magnetic Reynolds number, $S$ the coupling number, $\f: \Omega \to \R^3$ a source term and $\eps \u = \frac{1}{2} (\grad \u + \grad \u^\top )$. The system is completed with the boundary conditions
\begin{equation}\label{eq:boundarycond3d}
	\u=\mathbf{0}, \quad  \E \times \n = \mathbf{0}, \quad \B \cdot \n=0 \quad \text{ on } \del \Omega,
\end{equation}
where $\n$ is the unit outer normal vector.
The above formulation based on the electric and magnetic fields was first rigorously analyzed by Hu et al.\ \cite{Hu2020}.

In two dimensions, the electric field is a scalar field and hence denoted as $E$. The curl operators and cross products are interpreted depending on whether the arguments are scalar- or vector-valued corresponding to the definitions
\begin{equation}
	\scurl \B = \del_x B_2 - \del_y B_1, \qquad  \vcurl E = \begin{pmatrix}
		\del_y E\\
		-\del_x E
	\end{pmatrix}
\end{equation}
and 
\begin{equation}
	\u \times \B =u_1 B_2 - u_2 B_1, \qquad
	\B \times E =\begin{pmatrix}
		B_2 E \\
		-B_1 E
	\end{pmatrix}.
\end{equation}
Moreover, the boundary conditions for the electric field change to $E=0$ on $\del \Omega$ in two dimensions.

Other formulations include the current density $\mathbf{j}= \E + \u \times \B$ \cite{Hu2018} as an unknown or eliminate the electric field using equation \eqref{eq:MHD3d3}. In addition to the stationary case, we also consider the time-dependent version of \eqref{eq:MHD3d} where  the time-derivatives $\frac{\del\u}{\del t}$ and $\frac{\del\B}{\del t}$ are added to \eqref{eq:MHD3d1} and \eqref{eq:MHD3d4} respectively with suitable initial conditions $\u(\x,0)=\u_0(\x)$ and $\B(\x,0)=\B_0(\x) \, \forall \x \in \Omega$. Note that MHD models neglect displacement currents $\frac{\del\E}{\del t}$ \cite[Sec. 1.5]{Gerbeau2006}.

The main contribution of this work is to provide block preconditioners for linearizations of \eqref{eq:MHD3d} with good convergence even at high Reynolds and coupling numbers. The performance relies on the following three (novel) approaches:
\begin{itemize}
	\item[1.)] We consider a fluid-Reynolds-robust augmented Lagrangian preconditioner for an $\hdiv\times L^2$-discretization of the Navier--Stokes equations that relies on a specialized multigrid method.
	\item[2.)]  We introduce a new monolithic multigrid method for the electromagnetic block.
	\item[3.)]  We discover that using the outer Schur complement which eliminates the $(\u,p)$ block instead of the $(\E,\B)$ block has crucial advantages for ensuring robustness for high parameters.
\end{itemize}

Furthermore, we show that our preconditioners extend in a straightforward manner to the time-dependent version of \eqref{eq:MHD3d}. This has the substantial advantage that the choice of the time-stepping scheme is no longer restricted by the ability to solve the linear systems. In particular, it allows the use of fully implicit methods for high Reynolds numbers and coupling parameters.

An important point for discretizations is the enforcement of the magnetic Gauss' law $\div \B=0$ in the weak formulation, achieved in most cases by a non-physical Lagrange multiplier $r$ \cite{Schoetzau2004}. However, in general a Lagrange multiplier only enforces the divergence constraint in a weak sense, which can cause severe problems for the discretization and numerical simulations \cite{Brackbill1980, Dai1998}. For the $\B$-$\E$ formulation \eqref{eq:MHD3d} Hu et al.~\cite{Hu2020} show that both a Lagrange multiplier and an augmented Lagrangian term lead to a point-wise preservation of Gauss' law with appropriate choices of discrete spaces. In this work, we consider the latter approach by replacing \eqref{eq:MHD3d4} with
\begin{equation}
\frac{1}{\Rem}\grad \div\B + \vcurl \E = \mathbf{0}.
\end{equation}

The literature proposes numerous numerical schemes and preconditioning strategies for the numerical solution of the different formulations. The most common approach is based on block preconditioners in both the stationary \cite{Li2017, PhillipsPHD, Phillips2014, Wathen2017, Wathen2020} and time-dependent \cite{Chacn2008,Cyr2013, Phillips2016} cases. Here, the main challenges are to find suitable approximations of one or more Schur complements and robust linear solvers for the inner auxiliary problems. Phillips et al.~\cite{Phillips2016} simplify the Schur complement by the use of vector identities and approximate the remaining parts based on a spectral analysis. They report iteration counts for a stationary three-dimensional lid-driven cavity problem up to $\Re=\Rem=100$. A similar approach is used by Wathen and Greif in \cite{Wathen2020} where they construct an approximate inverse block preconditioner by sparsifying a derived formula for the exact inverse and drop low order terms. Here, results for Hartmann numbers $\Ha=\sqrt{S\Rem\Re}$ up to 1,000 are reported for stationary problems. Other approaches include fully-coupled geometric \cite{adler2020monolithic, Adler2016} and algebraic \cite{Shadid2010, Shadid2016} monolithic multigrid methods. In \cite{adler2020monolithic}, Adler et al.\ present results for a two-dimensional Hartmann problem for parameters up to $\Re = \Rem = 64$.

However, the performance of most of these preconditioners deteriorates significantly for high Reynolds and coupling numbers. To the best of our knowledge, a practical robust preconditioner for the stationary MHD equations has not yet been proposed. The common problem for high magnetic Reynolds numbers and coupling numbers for the stationary case is that all available Schur complement approximations become less accurate for Newton-type linearizations, causing the linear solver to fail to converge. Conversely, Picard-type linearizations can allow an exact computation of the Schur complement but fail to converge in the nonlinear iteration.

In this work, we consider two different linearizations. The first is the Picard iteration proposed by Hu et al.~\cite{Hu2020}. We compute an approximation to the outer Schur complement of the arising block system and introduce a robust linear solver for the different blocks. This scheme works well for small magnetic Reynolds numbers but the nonlinear iteration fails to converge for higher $\Rem$, as anticipated in the analysis of \cite{Hu2020}. The second is a full Newton linearization, which converges well for high Reynolds numbers and coupling numbers. However, our approximation of the Schur complement deteriorates slightly for high parameters.

Ma et al.~\cite{Ma2016} have developed Reynolds-robust preconditioners for the time-depen\-dent MHD equations that are based on norm-equivalent and field-of-values equivalent approaches. To the best of our knowledge, their strategy does not extend to the stationary case; in general, the time-dependent case offers crucial advantages for the development of robust solvers. For example, Ma et al.\ treat complicated terms like the hydrodynamic convection term $\u \cdot \grad \u$ explicitly in the time-stepping scheme, which can cause problems for convection-dominated problems and does not apply in the stationary case. The discretization of the time derivative causes mass matrices with a scaling of $1/\Delta t$, where $\Delta t$ denotes the time step size, to appear in the block matrix on the diagonal blocks. As we will see also in our numerical results for the time-dependent problems, these extra terms dominate the scheme for small $\Delta t$ and hence simplify the development of robust solvers.

Most applications are in the regime of high Reynolds and coupling numbers and hence it is of great interest to build robust solvers with respect to these parameters.
For liquid metals, the fluid Reynolds number $\Re$ tends to be much larger than $\Rem$. For example, the flow of liquid mercury is characterized by a ratio of $10^7$ between these two constants; typical values in aluminium electrolysis are $\Rem = 10^{-1}$ and $\Re = 10^5$ \cite{Gerbeau2006}. High magnetic Reynolds numbers occur on large length scales, as in geo- and astrophysics. The magnetic Reynolds number of the outer Earth's core is in the range of $10^3$ and of the sun is in the range of $10^6$ \cite{Davies2015}. Magnetic Reynolds numbers between $10^1-10^3$ have been used in several dynamo experiments that investigate planetary magnetic fields \cite{Molokov2007}. The coupling number $S$ is around $10^0$ for aluminium electrolysis \cite{Gerbeau2006} and Armero \& Simo \cite{Armero1996} define strong coupling for $S$ in the range of $10^2-10^9$.


The remainder of this work is outlined as follows.
In Section~\ref{sec:discretization}, we derive an augmented Lagrangian formulation for \eqref{eq:MHD3d} and describe the finite element discretization and linearization schemes. In Section~\ref{sec:derivationofblockpreconditioner}, we introduce block preconditioners for these schemes, present a calculation of the corresponding (approximate) Schur complements and describe robust linear multigrid solvers for the different blocks. Numerical examples and a detailed description of the algorithm are presented in Section~\ref{sec:numericalresults}.

\section{Formulation, linearization, and discretization} \label{sec:discretization}

\subsection{An augmented Lagrangian formulation}\label{sec:AnaugmentedLagrangianFormulation}
We modify \eqref{eq:MHD3d} by introducing two augmented Lagrangian terms: $-\gamma \grad \div \u$ for $\gamma>0$ is added to \eqref{eq:MHD3d1}, and (as previously discussed) $-1/\Rem\ \grad \div \B$ is added to \eqref{eq:MHD3d4}. Note that both terms leave the continuous solution of the problem unchanged. We use the first term to control the Schur complement of the fluid subsystem \cite{BenziOlshanskii,Farrell2020} and the second term to enforce the divergence constraint $\div \B=0$, as shown in \cite[Thm. 9]{Hu2020}.

Following these modifications, we consider the following system
\begin{subequations}
	\label{eq:MHDFinal}
	\begin{align}
		- \frac{2}{\mathrm{Re}} \div \mathbf{\varepsilon}( \u)  + \u \cdot \nabla \u - \gamma \grad \div \u  + \nabla p + S\, \B \times  (\E + \u \times \B) &= \f, \label{eq:MHDFinal1}\\
		\div \u&=0, \label{eq:MHDFinal2}\\
		\E + \u \times \B - \frac{1}{\Rem} \vcurl \B &= \mathbf{0}, \label{eq:MHDFinal3}\\
		-\frac{1}{\Rem} \grad \div \B + \vcurl \E &= \mathbf{0} \label{eq:MHDFinal4},
	\end{align}
\end{subequations}
subject to the boundary conditions \eqref{eq:boundarycond3d}.
For convenience, we consider homogeneous boundary conditions in this section but all the results extend in a straightforward manner to inhomogeneous boundary conditions. However, there are subtle technicalities for the implementation of the degrees of freedom in the finite element method in the inhomogeneous case, which are explained in detail in Section~\ref{sec:interp-boundary-data}.

The weak formulation of \eqref{eq:MHDFinal} seeks $\UU \coloneqq(\u,p,\E,\B)\in \ZZ \coloneqq \V\times Q \times \Rr \times \W$ with
\begin{equation}
	\V \coloneqq \Hozv, \quad Q\coloneqq L^2_0(\Omega),\quad \Rr \coloneqq \hzcurl, \quad  \W \coloneqq \hzdiv.
\end{equation}
In two dimensions,  the space for the electric field is scalar-valued and can be identified with $R\coloneqq \Hoz$.
The weak formulation is to find $\mathcal{U} \in \ZZ$ such that for all $\mathcal{V}\coloneqq(\v,q,\F,\C) \in \ZZ$ and $\mathcal{F} =  (\f, 0, \mathbf{0}, \mathbf{0})$ there holds
\begin{equation} \label{eq:weakform}
 \mathcal{R}(\mathcal{U}, \mathcal{V})	\coloneqq \NN(\mathcal{U}, \mathcal{V}) - (\mathcal{F},\mathcal{V}) = 0
\end{equation}
with
\begin{align}
	\begin{split}
		\NN(\mathcal{U},\mathcal{V}) & =\frac{2}{\Re}(\mathbf{\varepsilon}(\u), \mathbf{\varepsilon}(\v)) + (\u\cdot \grad \u, \v) + \gamma (\div \u, \div \v) \\
		& - (p,\div \v) +  S (\B \times \E, \v) + S (\B \times (\u \times \B), \v)  \label{eq:NN3} \\
		&- (\div \u, q) \\
		&+ (\E, \F) + (\u\times \B, \F) - \frac{1}{\Rem}(\B, \vcurl \F) \\
		&+ \frac{1}{\Rem} (\div \B, \div \C) + ( \vcurl \E, \C).
	\end{split}
\end{align}
All boundary integrals that result from integration by parts vanish because of the choice of the boundary conditions \eqref{eq:boundarycond3d}.

Note that $\W$ and $\Rr$ are chosen from the de Rham complex \cite{Arnold2018}
\begin{equation}\label{eq:contDC3d}
	\R \xrightarrow[]{\text{id}} \Hoz \xrightarrow[]{\mathrm{grad}} \hzcurl\xrightarrow[]{\vcurl} \hzdiv \xrightarrow[]{\div} L^2_0(\Omega) \xrightarrow[]{\text{null}} 0,
\end{equation}
which is exact for the simply connected domains we consider.
The corresponding complex in two dimensions is given by
\begin{equation}\label{eq:contDC2d}
	\R \xrightarrow[]{\text{id}} \mathrm{H}_0(\mathrm{curl}, \Omega) \xrightarrow[]{\vcurl} \hzdiv \xrightarrow[]{\div} L^2_0(\Omega) \xrightarrow[]{\text{null}} 0.
\end{equation}
This ensures that formulation \eqref{eq:MHDFinal} enforces the divergence constraint $\div \B =0$ and $\vcurl \E = \mathbf{0}$. To see this, we test \eqref{eq:weakform} with $\mathcal{V} = (\mathbf{0}, 0, \mathbf{0}, \vcurl \E)$ and conclude that $\vcurl \E =\mathbf{0}$. Here, $\mathcal{V}$ is a valid test function because the above exact sequence implies that $\vcurl(\Rr) = \W$. Similarly, testing with $\mathcal{V} = (\mathbf{0}, 0, \mathbf{0}, \B)$ results in $\div \B = 0$.

\subsection{Linearization: Newton and Picard}

The Newton linearization of \eqref{eq:weakform} for the initial guess $\mathcal{U}^n = (\u^n,p^n,\E^n,\B^n)$ is to find an update $\delta \mathcal{U}$ such that
\begin{align}
	\NN_\text{N}(\delta\mathcal{U},\mathcal{U}^n,\mathcal{V})&=\mathcal{R}(\mathcal{U}^n,\mathcal{V}) \quad \forall \ \mathcal{V} \in \ZZ,\\
	\mathcal{U}^{n+1} &= \mathcal{U}^n + \delta\mathcal{U},
\end{align}
with the weak form of the nonlinear residual $\mathcal{R}(\mathcal{U}^n, \mathcal{V})$ evaluated at $\mathcal{U}^n$ and

\begin{align}
	\label{eq:Newton}
	\begin{split}
		\NN_\text{N}(\delta\mathcal{U},\mathcal{U}^n,\mathcal{V}) &= \frac{2}{\Re}(\mathbf{\varepsilon}( \delta\u), \mathbf{\varepsilon}(\v)) + (\u^n\cdot \grad \delta\u, \v) + (\delta\u\cdot \grad \u^n, \v)  \\
		&+ \gamma (\div \delta\u, \div \v) - (\delta p,\div \v)\\
		& +  S (\B^n\times \delta \E, \v)  + S (\delta\B \times \E^n, \v)  \\
		& + S  (\B^n \times (\delta \u \times \B^n), \v) +
		S (\delta\B \times (\u^n \times \B^n), \v) \\
		&+ S (\B^n \times  (\u^n \times \delta\B), \v)\\
		&- (\div \delta\u, q) \\
		& + (\delta \E, \F) + (\u^n\times\delta\B, \F) + (\delta\u \times \B^n, \F) \\
		&- \frac{1}{\Rem} (\delta \B, \vcurl \F) \\
		&+ \frac{1}{\Rem} (\div \delta \B, \div \C) + (\vcurl \delta \E, \C).
	\end{split}
\end{align}

The bilinear form for the Picard iteration that we consider is given by
\begin{align}
	\label{eq:Picard}
	\begin{split}
		\NN_\text{P}(\delta\mathcal{U},\mathcal{U}^n,\mathcal{V}) = & \NN_{N} (\delta\mathcal{U},\mathcal{U}^n,\mathcal{V}) - S (\delta\B \times \E^n, \v) - S (\B^n \times  (\u^n \times \delta\B), \v)\\
		& - S (\delta\B \times (\u^n \times \B^n), \v) - (\u^n\times\delta\B, \F).
	\end{split}
\end{align}
Note that in contrast to definition of the Picard iteration in \cite{Hu2020}, we do not scale the term $(\vcurl \delta \E, \C)$ with $S/\Rem$ and consider the full Newton linearization of the advection term $(\u\cdot \grad) \u$.
The advantage of the Picard linearization \eqref{eq:Picard} in comparison to the Newton linearization \eqref{eq:Newton} is that it allows an exact Schur complement computation in two dimensions and converges well for high $\Re$. However, its major disadvantage is the failure of nonlinear convergence for high $\Rem$.

\subsection{Discretization}

For a finite element discretization, we seek $\mathcal{U}_h:=(\u_h, p_h,$ $\E_h, \B_h) \in \ZZ_h \coloneqq \mathbf{V}_h \times Q_h \times \mathbf{R}_h\times \W_h$ such that
\begin{equation} \label{eq:weakformdiscr}
	\NN(\mathcal{U}_h, \mathcal{V}_h)=  (\mathcal{F},\mathcal{V}_h) \quad \forall\, \mathcal{V}_h \in \ZZ_h.
\end{equation}
We choose Raviart--Thomas elements of degree $k$ $\mathbb{RT}_k$ \cite{Raviart1977} for $\W_h$, \Ned elements of first kind $\mathbb{NED}1_k$ \cite{Nedelec1980} for $\Rr_h$ in 3D and continuous Lagrange elements $\mathbb{CG}_k$ for $R_h$ in 2D. Note that these elements belong to a discrete subcomplex of \eqref{eq:contDC3d}
\begin{equation}\label{eq:deRhamdiscret3d}
	\mathbb{CG}_k  \xrightarrow[]{\mathrm{grad}} \mathbb{NED}1_k \xrightarrow[]{\vcurl} \mathbb{RT}_k \xrightarrow[]{\div} \mathbb{DG}_{k-1} \xrightarrow[]{\text{null}} 0,
\end{equation}
and of \eqref{eq:contDC2d}
\begin{equation}\label{eq:deRhamdiscret2d}
	\mathbb{CG}_k \xrightarrow[]{\vcurl} \mathbb{RT}_k \xrightarrow[]{\div} \mathbb{DG}_{k-1} \xrightarrow[]{\text{null}} 0.
\end{equation} This implies that we enforce $\div \B_h=0$ and $\vcurl \E_h=\mathbf{0}$ pointwise with the same proof as for the continuous case. These identities also hold for inhomogeneous boundary conditions, since the interpolation operator $\II^h_{\W_h}$ into the Raviart--Thomas space satisfies for all divergence-free $\B \in \mathbf{H}_0(\mathrm{div}, \Omega)$~\cite[Prop.~2.5.2]{Boffi2013}
\begin{equation}\label{eq:InterpolationPreserveRT}
	\div(\II^h_{\W_h} \B) = 0.
\end{equation}

Moreover, following \cite{douglas1976}, we add the following stabilization term to address the problem that the Galerkin discretization of advection-dominated problems problems can be oscillatory \cite{Elman2014}
\begin{equation}\label{eq:stabBurman}
 \sum_{K\in \MM_h} \frac{1}{2} \int_{\del K} \mu \, h_{\del K}^2 \llbracket \grad \u_h \rrbracket : \llbracket \grad \v_h \rrbracket \ \mathrm{d} s.
\end{equation}
Here, $\llbracket \grad \u_h \rrbracket$ denotes the jump of the gradient, $h_{\del K}$ is a function giving the facet size, and $\mu$ is a free parameter that is chosen according to \cite{burman_edge_2006}.

Note that a fully robust discretization should also include a stabilization term for the magnetic field $\B$ in the case of dominating magnetic advection. The literature does not propose many stabilization types for this problem. The most promising work by Wu and Xu \cite{Wu2020} uses the so-called SAFE-scheme for stabilization which is based on an exponential fitting approach. While the original SAFE-scheme is only a first order method, it can be extended to higher order as shown in \cite{wu2020unisolvence}. We aim to include this stabilization in future work.

For the hydrodynamic part, we consider the $\hdiv\times L^2$-conforming element pair $\mathbb{BDM}_k\times \mathbb{DG}_{k-1}$ with the Brezzi-Douglas-Marini element $\mathbb{BDM}_k$ of order $k$ \cite{Brezzi1985, Ndlec1986}. This discretization ensures that $\div \u_h =0$ holds pointwise since $\div \mathbf{V}_h \subset Q_h$. Additionally, it exhibits pressure robustness, i.e.\ the error estimates do not degrade for high Reynolds numbers \cite{john2017}.

Since the discretization is nonconforming, we must consider a discontinuous Galer\-kin formulation of the hydrodynamic advection and diffusion terms \cite[section 7]{GauerLinke2019}. We denote by $\FF_h = \FF_h^i \cup \FF_h^\partial$ all facets of the triangulation, which consists of the interior facets $\FF_h^i $ and the Dirichlet boundary facets $\FF_h^\partial$. We assign to each $F \in \mathcal{F}_h$ its diameter $h_F$ and unit normal vector $\mathbf{n}_F$. The jump and average operators across a facet are denoted by $\llbracket\cdot \rrbracket$ and $\ldblbrace\cdot \rdblbrace$ respectively and are defined as $\llbracket \Phi \rrbracket=\Phi^+ - \Phi^-$ and $\ldblbrace\Phi \rdblbrace=\frac{1}{2}(\Phi^+ + \Phi^-)$. The penalization parameter is chosen as $\sigma = 10 k^2$. Inhomogeneous boundary data are described by $\mathbf{g}_D$. We then add the following bilinear forms to \eqref{eq:weakformdiscr}:
\begingroup
\begin{align}
	\begin{split}
		a_h^{DG}(\u_h, \v_h)= &-\frac{2}{\Re}\sum_{F \in \FF_h} \int_F \ldblbrace\varepsilon( \u_h) \rdblbrace \mathbf{n}_F \cdot \llbracket \v_h \rrbracket \,\mathrm{d}s\\
		&-\frac{2}{\Re}\sum_{F \in \FF_h} \int_F  \llbracket \u_h \rrbracket \cdot \ldblbrace\varepsilon(\v_h) \rdblbrace \mathbf{n}_F  \,\mathrm{d}s\\
		&+\frac{1}{\Re}\sum_{F \in \FF_h} \frac{\sigma}{h_F} \int_F \llbracket\u_h \rrbracket \cdot \llbracket \v_h \rrbracket \,\mathrm{d}s  \\
		& - \frac{1}{\Re}\sum_{F \in \FF_h^\partial} \frac{\sigma}{h_F} \int_F \mathbf{g}_D \cdot \v_h \,\mathrm{d}s\ +  \frac{2}{\Re}\sum_{F \in \FF_h^\partial} \int_F \mathbf{g}_D \cdot \varepsilon( \v_h) \mathbf{n}_F \,\mathrm{d}s,
	\end{split}
	\\
	\begin{split}
		c_h^{DG}(\u_h, \v_h)= &\ \ \ \ \frac{1}{2}\sum_{F \in \FF^i_h}\int_F \llbracket (\u_h \cdot \mathbf{n}_F + |\u_h\cdot \n_F|)\u_h \rrbracket \cdot
		\llbracket \v_h \rrbracket \,\mathrm{d}s  \\
		& +\frac{1}{2} \sum_{F \in \FF_h^\partial} \int_F (\u_h \cdot \mathbf{n}_F + |\u_h\cdot \n_F|) \u_h \cdot \v_h \,\mathrm{d}s \\
		& +\frac{1}{2} \sum_{F \in \FF_h^\partial} \int_F (\u_h \cdot \mathbf{n}_F - |\u_h\cdot \n_F|) \mathbf{g}_D \cdot \v_h \,\mathrm{d}s.
	\end{split}
\end{align}
\endgroup


Hu et al.\ prove in \cite[Theorem 4]{Hu2020} that \eqref{eq:weakformdiscr} is well-posed and has at least one solution. The solution is unique for suitable source and boundary data.
While the well-posedness and convergence of the Newton iteration remains an open problem, Hu et al.\ prove that the Picard iteration converges to the unique solution of \eqref{eq:weakformdiscr} if both $\Re^2 \|f\|_{-1}$ and $\Re\Rem^\frac{3}{2} \|f\|_{-1}$ are small enough.

For the Newton linearization \eqref{eq:Newton}, we must solve the following linear system at each step:
\begin{equation}
	\label{eq:matrix_upBE}
	\begin{bmatrix}
		\mathcal{F} +\DD & \BB^\top &  \JJ & \tilde{\JJ} + \tilde{\DD}_1 + \tilde{\DD}_2 \\
		\BB & \mathbf{0} & \mathbf{0} & \mathbf{0} \\
		\GG& \mathbf{0}& \MM_\E & \tilde{\GG} -\frac{1}{\Rem} \AA \\
		\mathbf{0} & \mathbf{0} & \AA^\top  & \CC
	\end{bmatrix}
	\begin{bmatrix}
		x_\u \\ x_p \\ x_\E \\ x_\B
	\end{bmatrix} =
	\begin{bmatrix}
		\mathcal{R}_\u\\ \mathcal{R}_p\\ \mathcal{R}_\E\\ \mathcal{R}_\B
	\end{bmatrix},
\end{equation}
where $x_\u$, $x_p$, $x_\E$ and $x_\B$ are the coefficients of the discretized Newton corrections and $\mathcal{R}_\u$, $\mathcal{R}_p$, $\mathcal{R}_\E$ and $\mathcal{R}_\B$ the corresponding nonlinear residuals.
The correspondence between the discrete and continuous operators is illustrated in Table \ref{tab:Operators}. We have chosen the notation that operators that include a tilde are omitted in the Picard linearization $\eqref{eq:Picard}$. Moreover, we introduce $\eta\in\{0,1\}$ to distinguish between the stationary ($\eta=0$) and transient ($\eta=1$) cases.

For the time-dependent equations, we concentrate here on the implicit Euler method, but the following computations are straightforward to adapt to other implicit multi-step methods. We use the same finite element discretization as in the stationary case. Note that in the transient case, the equation
\begin{equation}
	\del_t \B + \vcurl \E = \mathbf{0}
\end{equation}
immediately implies $\div \B=0$ if the initial condition satisfies $\div \B_0=0$, and this remains true on the discrete level up to solver tolerances; see \cite[Theorem 1]{Hu20162} for a proof for implicit Euler which can be extended to other multi-step methods in a straightforward manner, provided all starting values are divergence-free.
Hence, the augmented Lagrangian term $-	\frac{1}{\Rem}\grad \div \B$ is no longer necessary to enforce the divergence constraint and could therefore be omitted. Nevertheless, we retain it in our scheme since we employ the identity
\begin{equation}\label{eq:laplaceidentity}
	\frac{1}{\Rem} \vcurl \vcurl \u - 	\frac{1}{\Rem} \grad \div \u = - 	\frac{1}{\Rem} \Delta \u
\end{equation}
in our derivation of Schur complement approximations below.

\begin{table}[htb!]
  \centering
  \resizebox{\textwidth}{!}{
			\begin{tabular}{c|c|c}
				\toprule
				\textbf{Discrete} & \textbf{Continuous} & \textbf{Weak form}\\
				\midrule
				$\mathcal{F} \u$ & $\frac{\eta}{\Delta t} \u-\frac{2}{\Re} \div \varepsilon(\u) + \u^n\cdot \grad \u $ & $\frac{\eta}{\Delta t} (\u, \v)+\frac{2}{\Re}(\varepsilon( \u), \varepsilon( \v)) + (\u^n\cdot \grad \u , \v)$  \\
				& $+ \u \cdot \grad \u^n-\gamma\grad \div \u$ &  $+(\u\cdot \grad \u^n, \v) + \gamma(\div \u, \div \v)  $  \\
				$\DD \u$& $S \B^n\times(\u\times\B^n)$ & $ S (\B^n\times (\u \times \B^n), \v)$ \\
				$\JJ \E$ & $ S \B^n\times\E$ & $ S (\B^n\times\E, \v)$   \\
				$\tilde{\JJ} \B$ & $ S \B\times\E^n$ & $ S (\B\times\E^n, \v)$   \\
				$\tilde{\DD}_1\B$ & $ S \B \times (\u^n\times\B^n)$ & $ S (\B \times (\u^n\times\B^n),\v)$  \\
				$\tilde{\DD}_2\B$ & $ S \B^n \times (\u^n\times\B)$ & $ S (\B^n \times (\u^n\times\B),\v)$  \\
				$\MM_\E \E$ & $\E$ & $(\E,\F)$  \\
				$\GG \u$ & $\u \times \B^n$ & $(\u\times\B^n, \F)$  \\
				$\tilde{\GG} \B$ & $\u^n \times \B$ & $(\u^n\times\B, \F)$  \\
				$\AA \B$ & $ \vcurl \B$ & $ (\B, \vcurl \F)$ \\
				$\CC \B$ & $\frac{\eta}{\Delta t} \B-\frac{1}{\Rem}\grad \div \B$ & $\frac{\eta}{\Delta t}(\B, \C) + \frac{1}{\Rem}(\div\B,\div\C)$  \\
				$\AA^\top \E$ & $\vcurl \E$ & $(\vcurl \E, \C)$ \\
				$\BB^\top p$ & $\grad p$ & $-(p, \div \v)$ \\
				$\BB \u$ & $-\div \u$ & $-(\div \u, q)$  \\
				\bottomrule
		\end{tabular}}
		\caption{Overview of operators. Terms that include a tilde are dropped in the Picard iteration. The stationary and transient cases are distinguished by $\eta\in\{0,1\}$. }
		\label{tab:Operators}
\end{table}

\section{Derivation of block preconditioners}\label{sec:derivationofblockpreconditioner}
We now consider block preconditioners for \eqref{eq:matrix_upBE}. The inverse of a 2x2 block matrix can factorized as \cite{Benzi2005,Elman2014} 
\begin{equation}\label{eq:blockfactor}
	\begin{pmatrix} \MM & \mathcal{K} \\ \mathcal{L} & \mathcal{N} \end{pmatrix}^{-1}
	=
	\begin{pmatrix}
		\II &-\MM^{-1}\mathcal{K} \\
		\mathbf{0} &\II
	\end{pmatrix}
	\begin{pmatrix}
		\MM^{-1} &\mathbf{0} \\
		\mathbf{0} &\mathcal{S}^{-1}
	\end{pmatrix}
	\begin{pmatrix}
		\II &\mathbf{0} \\
		-\mathcal{L}\MM^{-1} &\II
	\end{pmatrix}
\end{equation}
provided the top-left block $\MM$ and the Schur complement $\mathcal{S}=N - \mathcal{L}\MM^{-1} \mathcal{K}$ are invertible. Since the Schur complement is usually a dense matrix, the main task is find a suitable approximation $\tilde{\mathcal{S}}$ for the Schur complement $\mathcal{S}$ as well as efficient solvers for $\MM$ and $\tilde{\mathcal{S}}$.

In Sections \ref{sec:OuterSchurEB} and \ref{sec:OuterSchurup} we derive approximations of the Schur complements for two different block elimination strategies.  We briefly introduce the theory of parameter-robust multigrid relaxation in Section \ref{sec:robustssc}, and then describe the multigrid methods that we use to solve the top-left block and the Schur complement approximations in Sections \ref{sec:solverforschurcomp} and \ref{sec:solverformagneticblock}.

Both block preconditioners we consider gather the variables as $(\E_h, \B_h)$ and $(\u_h, p_h)$. They differ in the order of block elimination: the first takes the Schur complement that eliminates (inverts) the $(\E_h, \B_h)$ block, while the second takes the Schur complement that eliminates the $(\u_h, p_h)$ block. The first choice appears several times in the literature \cite{Li2017,Phillips2016}, while it seems that the second choice has not yet been investigated. As we will see, for small $\Rem$ and $S$ both preconditioners perform similarly, while for more difficult parameter regimes the second choice notably outperforms the first. We therefore recommend the second strategy and mainly report numerical results for this choice. Nevertheless, we also investigate the first option, both for comparison and because it allows a much more detailed description of the Schur complement. In two dimensions it even allows an exact computation of the Schur complement. The two strategies are compared in Section \ref{sec:stationaryliddrivencavityproblemin3d} below.

\subsection{Outer Schur complement eliminating the $(\E_h, \B_h)$ block}\label{sec:OuterSchurEB}
Reordering \eqref{eq:matrix_upBE} for convenience, we consider
\begin{equation}
	\label{eq:matrix_EBup}
	\left[
	\begin{array}{cc|cc}
		\MM_\E & \tilde{\GG} - \frac{1}{\Rem} \AA & \GG & \mathbf{0}\\
		\AA^\top & \CC & \mathbf{0} & \mathbf{0} \\
		\hline \rule{0pt}{1.0\normalbaselineskip}
		\JJ & \tilde{\JJ} + \tilde{\DD}_1 + \tilde{\DD}_2 & \mathcal{F} + \DD & \BB^\top\\
		\mathbf{0} & \mathbf{0} & \BB & \mathbf{0}
	\end{array}
	\right]
	\begin{bmatrix}
		x_\E \\ x_\B \\ x_\u \\ x_p
	\end{bmatrix} =
	\begin{bmatrix}
		\mathcal{R}_\E\\ \mathcal{R}_\B\\ \mathcal{R}_\u\\ \mathcal{R}_p
	\end{bmatrix}.
\end{equation}
In the following, we refer to the Schur complement of the 4x4 matrix as the outer Schur complement, while we call the Schur complements of the resulting 2x2 blocks inner Schur complements.
The outer Schur complement eliminating the $(\E_h, \B_h)$ block is given by
\begin{equation}
	\label{eq:outerSchurCompNewton}
	\mathcal{S}^{(\E, \B)} =
	\begin{bmatrix}
		\FF+\DD & \BB^\top \\
		\BB & \mathbf{0}
	\end{bmatrix}
	-
	\begin{bmatrix}
		\JJ & \tilde{\JJ} + \tilde{\DD}_1 + \tilde{\DD}_2\\
		\zerov & \zerov
	\end{bmatrix}
	\begin{bmatrix}
		\MM_\E& \tilde{\GG} - \frac{1}{\Rem} \AA \\
		\AA^\top & \CC
	\end{bmatrix}^{-1}
	\begin{bmatrix}
		\GG & \zerov\\
		\zerov & \zerov
	\end{bmatrix}.
\end{equation}

We simplify $\mathcal{S}^{(\E, \B)}$ by applying the identity \eqref{eq:blockfactor}
to the top-left electromagnetic block
\begin{equation}
	\label{eq:M}
	\MM =
	\begin{bmatrix}
		\MM_\E& \tilde{\GG} - \frac{1}{\Rem} \AA \\
		\AA^\top & \CC
	\end{bmatrix}.
\end{equation}
This results in
\begin{equation}
	\mathcal{S}^{(\E, \B)}=
	\begin{bmatrix}
		\FF + \DD - \JJ\MM^{-1}_{1,1} \GG - (\tilde{\JJ} + \tilde{\DD}_1 + \tilde{\DD}_2)\MM^{-1}_{2,1}\GG& \BB^\top \\
		\BB & \zerov
	\end{bmatrix}
\end{equation}
with
\begin{equation}
	\MM^{-1}_{1,1} = \MM_\E^{-1} + \MM_\E^{-1}\left(\tilde{\GG}-\frac{1}{\Rem} \AA\right)\left(\CC - \AA^\top \MM_\E^{-1}\left(\tilde{\GG}-\frac{1}{\Rem} \AA\right)\right)^{-1}\AA^\top\MM_\E^{-1}
\end{equation}
and
\begin{equation}
	\MM^{-1}_{2,1} = -\left(\CC - \AA^\top \MM_\E^{-1}\left(\tilde{\GG}-\frac{1}{\Rem} \AA\right)\right)^{-1}\AA^\top\MM_\E^{-1}.
\end{equation}

We precondition $\mathcal{S}^{(\E, \B)}$ for both linearizations in the stationary case by 
\begin{equation}\label{eq:SchurStat}
	\tilde{\mathcal{S}}^{(\E, \B)} =
	\begin{bmatrix}
		\FF+\DD & \BB^\top \\
		\BB & \mathbf{0}
	\end{bmatrix},
\end{equation} 
and in the transient case by
\begin{equation}\label{eq:SchurAlpha}
	\tilde{\mathcal{S}}^{(\E, \B)}_{\alpha}:=
	\begin{bmatrix}
		\FF+\alpha\DD & \BB^\top \\
		\BB & \mathbf{0}
	\end{bmatrix},
	\quad \quad \alpha = \frac{\Delta t}{\Delta t + \Rem h^2 + \delta \Rem h \|\u^n\|_{L^2}\Delta t }.
\end{equation}

In the following, we motivate this choice of preconditioners and emphasize the cases in which these Schur complement approximations are exact. Therefore, we mainly follow \cite{Phillips2016}, but adapt the computations for our formulation which includes the electric field $\E$ instead of a Lagrange multiplier $r$.

For the simplification of the outer Schur complement $\mathcal{S}^{(\E, \B)}$ we must find approximations for
\begin{equation}\label{eq:SchurcompExtraterms}
\mathcal{K}_1 \coloneqq	\DD - \JJ\MM^{-1}_{1,1} \GG \quad \text{  and  } \quad  \mathcal{K}_2 \coloneqq -(\tilde{\JJ} + \tilde{\DD}_1 + \tilde{\DD}_2)\MM^{-1}_{2,1}\GG.
\end{equation}
Note that the first summand of $\JJ\MM^{-1}_{1,1} \GG $ is $\JJ \MM_\E^{-1} \GG$ which equals $\DD$. Hence, $\mathcal{K}_1$ simplifies to the second summand of  $\JJ \MM_\E^{-1} \GG$, i.e., 
\begin{equation}\label{eq:Schurcomppart1}
	\mathcal{K}_1 = -\JJ \MM_\E^{-1}\left(\tilde{\GG}-\frac{1}{\Rem} \AA\right)\left(\CC - \AA^\top \MM_\E^{-1}\left(\tilde{\GG}-\frac{1}{\Rem} \AA\right)\right)^{-1}\AA^\top\MM_\E^{-1}\GG
\end{equation}
which corresponds on a continuous level to
\begin{equation}\label{eq:Schurcomppart1cont}
	\scalemath{0.945}{
		-S\, \B^n\times \left( \left(\delta\,\u^n \times \cdot - \frac{1}{\Rem} \vcurl\right)\left(\frac{\eta}{\Delta t} I-\frac{1}{\Rem} \Delta - \delta\,\vcurl(\u^n \times \cdot)\right)^{-1}  \vcurl(\u \times \B^n)\right)},
\end{equation}
where  $\cdot$ denotes a placeholder for the input of the corresponding operators.
Moreover, we have used $\delta \in \{0,1\}$ to distinguish between the Picard ($\delta=0$) and Newton ($\delta=1$) linearizations. In the discrete counterpart \eqref{eq:Schurcomppart1}, the matrix arising in the Picard iteration is made by dropping all terms with a tilde.

The continuous expression for $\mathcal{K}_2$ is given by
\begin{align}\label{eq:Schurcomppart2cont}
	\begin{split}
		\delta\,S\,(\cdot \times \E^n + \cdot & \times(\u^n \times \B^n) + \B^n\times(\u^n\times \cdot)) \\
		& \left(\frac{\eta}{\Delta t}{I}-\frac{1}{\Rem}
		\Delta - \delta \vcurl(\u^n \times \cdot)\right)^{-1} \vcurl(\u \times \B^n).
	\end{split}
\end{align}
Note that $\mathcal{K}_2$ vanishes for the Picard iteration.

\subsubsection{The two-dimensional case}
For the Picard linearization, expression \eqref{eq:Schurcomppart1} simplifies to $\DD$ in the stationary case. This follows immediately from the two-dimen\-sional analogue of \eqref{eq:laplaceidentity} and the identity \cite{Phillips2014}
\begin{equation}\label{eq:curlidentity2d}
	\scurl (-\Delta)^{-1} \vcurl \varphi = \varphi
\end{equation}
which implies for our structure-preserving discretization that
\begin{equation}
	\AA (\CC + \AA^\top \MM_\E^{-1}\AA)^{-1}\AA^\top = \MM_\E.
\end{equation}
That means in the two-dimensional stationary case the outer Schur complement for the Picard iteration is exactly given by
\eqref{eq:SchurStat},
i.e.,~the Navier--Stokes block with the linearized Lorentz force.

In the transient case, the Schur complement for the Picard linearization can no longer be calculated exactly. The behavior of the Schur complement now depends on which of the terms $\frac{1}{\Delta t} I$ and $\frac{1}{\Rem} \Delta$ dominates in \eqref{eq:Schurcomppart1cont}. If $\frac{1}{\Delta t}$ is small in comparison to $\frac{1}{\Rem h^2}$, a good approximation of \eqref{eq:Schurcomppart1cont} is given, as in the stationary case, by $\tilde{\mathcal{S}}^{(\E, \B)}$. If $\frac{1}{\Delta t}$ dominates over $\frac{1}{\Rem h^2}$, \eqref{eq:Schurcomppart1cont}  is approximately given by
\begin{equation}
	S\, \B^n\times \left( \frac{1}{\Rem} \vcurl\, \left(\frac{1}{\Delta t} I \right)^{-1}  \vcurl(\u \times \B^n)\right).
\end{equation}
Hence, its magnitude can be approximated by $\frac{S \|\B^n\|^2\Delta t}{\Rem\, h^2}\ll1$ for moderate coupling numbers and therefore we neglect this term by using the approximation
	$\begin{bmatrix}
		\FF & \BB^T \\
		\BB & \mathbf{0}
	\end{bmatrix} $ for the Schur complement in this case.

 To also include the intermediate regime, we use the approximation of Phillips et al.\ \cite{Phillips2016} who suggest to use
\eqref{eq:SchurAlpha}.
The expression for $\alpha$ interpolates between the above mentioned dominating cases, since $\alpha \approx 0$ if $\frac{1}{\Delta t}\gg \frac{1}{\Rem h^2}$  and $\alpha \approx 1$ if  $\frac{1}{\Delta t}\ll \frac{1}{\Rem h^2}$.


For the Newton linearization, the simplification of $\mathcal{S}^{(\E, \B)}$ is not straightforward, but our numerical tests suggest that $\tilde{\mathcal{S}}^{(\E, \B)}$ and $	\tilde{\mathcal{S}}^{(\E, \B)}_{\alpha}$ are acceptable preconditioners for $\mathcal{S}^{(\E, \B)}$ in the stationary and transient cases, deteriorating only for high $S$ and $\Rem$.

This can be explained by the fact that for small $\Rem$ or $\Delta t$ the terms $\frac{1}{\Rem}\vcurl$ and $\frac{\eta}{\Delta t}{I}-\frac{1}{\Rem} \Delta$ dominate over $\delta \u^n\times \cdot$ and $\delta \vcurl(\u^n\times \cdot)$ in \eqref{eq:Schurcomppart1cont}. Remember that the terms that include a $\delta$ do not appear in the Picard iteration and were hence neglected in the previous derivation for the Picard iteration. Moreover, the term $\mathcal{K}_2$ is not included in our preconditioner for the Newton scheme which should deteriorate the performance for large $S$.


\subsubsection{The three dimensional case}\label{sec:derivationofblockpreconditioner3d}
The main difficulty in three dimensions is that the identity \eqref{eq:curlidentity2d} no longer holds. Therefore, $\tilde{\mathcal{S}}^{(\E, \B)}$ is not the exact outer Schur complement for the Picard linearization in the stationary case. In \cite{Phillips2016} the same approximation from the two dimensional case is used in three dimensions. Based on the argument for the two-dimensional case in the previous subsection, we expect this approximation to work well when the term $\Delta t$ dominates and to deteriorate in the other cases, especially in the stationary case.
The three dimensional performance of this preconditioner could be substantially improved with a better approximation of $\vcurl \Delta^{-1} \vcurl$ than a scaled identity.

We briefly comment on the main part of the outer Schur complement in the stationary case, given by
\begin{equation}\label{eq:SchurComp3dMain}
	S \B^n \times \left[\vcurl \Delta^{-1} \vcurl(\u \times \B^n) \right].
\end{equation}
As shown in \cite[Chapter 4]{PhillipsPHD} one can rewrite $\vcurl \Delta^{-1} \vcurl$ as $I -\nabla \Delta_r^{-1} \nabla \cdot$, where $\Delta_r$ denotes a scalar Laplacian. These two representations show that the operator is the identity on divergence-free functions and maps curl-free functions to zero. Hence, this operator corresponds to the orthogonal $L^2$-projection of a vector field onto its divergence-free part, which we denote by $\mathbb{P}$. Thus, the weak form of \eqref{eq:SchurComp3dMain} is given by
\begin{equation}\label{eq:SchurComp3dMainWeakForm}
	S(\mathbb{P}(\u \times \B^n), \mathbb{P}(\v \times \B^n)).
\end{equation}
 The key challenge is then to find a sparse approximation of \eqref{eq:SchurComp3dMainWeakForm}. We do not further address this challenge here and focus instead on the outer Schur complement that eliminates the $(\u_h,p_h)$ block.


%

\subsection{Outer Schur complement eliminating the $(\u_h, p_h)$ block}\label{sec:OuterSchurup}
The outer Schur complement eliminating the $(\u_h, p_h)$ is given by
\begin{equation}
	\label{eq:outerSchurCompNewton_BE}
	\mathcal{S}^{(\u, p)} =
	\begin{bmatrix}
		\MM_\E& \tilde{\GG} - \frac{1}{\Rem} \AA \\
		\AA^\top & \CC
	\end{bmatrix}
	-
	\begin{bmatrix}
		\GG & \zerov\\
		\zerov & \zerov
	\end{bmatrix}
	\begin{bmatrix}
		\FF+\DD & \BB^\top \\
        \BB & \mathbf{0}
	\end{bmatrix}^{-1}
	\begin{bmatrix}
		\JJ & \tilde{\JJ} + \tilde{\DD}_1 + \tilde{\DD}_2\\
		\zerov & \zerov
	\end{bmatrix}.
\end{equation}
Formula \eqref{eq:blockfactor} implies that the outer Schur complement for the Newton iteration is given by
\begin{equation}
	\mathcal{S}^{(\u, p)}=
	\begin{bmatrix}
		\MM_\E - \GG\NN^{-1}_{1,1} \JJ & \tilde{\GG} - \frac{1}{\Rem}\AA -\GG\NN^{-1}_{1,1}(\tilde{\JJ} + \tilde{\DD}_1 + \tilde{\DD}_2) \\
		\AA^\top &  \CC
	\end{bmatrix},
\end{equation}
where
\begin{equation}
    \NN^{-1}_{1,1} = (\FF+\DD)^{-1} - (\FF+\DD)^{-1}\BB^\top(-\BB(\FF+\DD)^{-1}\BB^\top)^{-1}\BB(\FF+\DD)^{-1}.
\end{equation}
For this strategy, further simplifications of the Picard or Newton linearizations are not straightforward. Our numerical results in the next section show that 
\begin{equation}
	\tilde{\mathcal{S}}^{(\u, p)} =
	\begin{bmatrix}
		\MM_\E& \tilde{\GG} - \frac{1}{\Rem} \AA \\
		\AA^\top & \CC
	\end{bmatrix}
\end{equation}
 works very well as a preconditioner for both schemes. Indeed, in contrast to the previous order of elimination, this approximation works qualitatively the same in two and three dimensions. We expect the approximation to deteriorate in the stationary case for very high $\Rem$, since the missing term $-\GG\NN^{-1}_{1,1}(\tilde{\JJ}+\tilde{\DD}_1 + \tilde{\DD}_2)$ in the Schur complement approximation gains more influence in comparison to $- \frac{1}{\Rem}\AA$. We also observe this numerically in the next section. 
 
 However, good approximation of the outer Schur complement is maintained for high coupling numbers $S$, which will clearly be seen in our numerical results in Section \ref{sec:stationaryliddrivencavityproblemin3d} below. This behavior is perhaps explained by the fact that $\NN^{-1}_{1,1}$ also includes a factor $S$ in the inverse of $(\FF+\DD)$, which cancels the factor of $S$ in the matrices $\JJ, \tilde{\JJ}, \tilde{\DD_1}$ and $\tilde{\DD_2}$.

To use these block preconditioners in practice, we must develop robust preconditioners for the electromagnetic and hydrodynamic subsystems.

\subsection{Parameter-robust relaxation}\label{sec:robustssc}

The equations in the hydrodynamic and electromagnetic blocks become difficult to solve in the parameter regimes of interest at high Reynolds and coupling numbers both due to the non-symmetric linearised advection and Lorentz force terms, and the addition of the symmetric positive semi-definite (SPSD) augmented Lagrangian terms. 
Standard multigrid methods are known to perform poorly for these kinds of problems. 

The key components for a robust multigrid method  for the SPSD augmented Lagrangian terms are a parameter-robust relaxation method, that efficiently damps error modes in the kernel of the singular operators, and a kernel-preserving prolongation operator, as revealed in the seminal work of Sch\"oberl~\cite{schoberl1999b}. The non-symmetric terms are more troublesome, but numerical results have shown \cite{Farrell2020,Farrell2018} that subspace correction methods can still perform well for the Navier--Stokes equations at high Reynolds numbers.

A recent summary of the theory of robust relaxation methods can be found in \cite{farrell2019pcpatch}. Briefly, we consider the multigrid relaxation methods in the framework of subspace correction methods \cite{xu1992}. These decompose a (finite-dimensional) trial space $V$ as
\begin{equation}\label{eq:decomp}
	V = \sum_{i} V_i,
\end{equation}
where the sum is not necessarily direct.
The parallel subspace correction method applied to a linear variational problem $a(u,v)=(f,v) \quad \forall v \in V$ computes for an initial guess $u^k$ an correction $\delta u_i$ to the error $e=u-u^k$ in each subspace $V_i$ by solving
\begin{equation}
    a(\delta u_i, v_i) = (f, v_i) - a(u^k, v_i) \text{ for all } v_i \in V_i,
\end{equation}
and sets $u^{k+1} = u^k + \sum_i w_i \delta u_i$ for damping parameters $w_i$.
A rigorous statement regarding the properties the decomposition \eqref{eq:decomp} and the considered bilinear form $a$ must fulfil to yield a robust relaxation method can be in found in \cite[Theorem 4.1]{schoberl1999b}.
A key property is that the kernel $\mathcal{N}$ of the SPSD terms is decomposed over the subspaces, i.e.,
\begin{equation}\label{eq:kerneldecomp}
	\mathcal{N} = \sum_i (V_i \cap \mathcal{N}).
\end{equation}
This property means that it must be possible to write any kernel function as the sum of kernel functions in the subspaces $V_i$. This implies that the subspaces $V_i$ must
be at least rich enough to support nonzero kernel functions. The choice of the space decomposition \eqref{eq:decomp} is often made by consideration of the
discrete Hilbert complexes underpinning the discretization.\\

\subsection{Solver for the hydrodynamic block}\label{sec:solverforschurcomp}


In order to implement the block factorization \eqref{eq:blockfactor} as the outer preconditioner, we need a solver for the Navier--Stokes subsystem. To do this,
we will apply ideas of parameter-robust multigrid relaxation described in Section \ref{sec:robustssc}, albeit without a theoretical guarantee of success. The
variational statement of the PDE we wish to solve is
	\begin{alignat}{2}
		\frac{2}{\Re}(\mathbf{\varepsilon}(\u), \mathbf{\varepsilon}(\v)) + (\u^n\cdot \grad \u, \v) +  (\u\cdot \grad \u^n, \v) + \gamma (\div \u, \div \v)& & \nonumber\\
	    + S (\B^n \times (\u \times \B^n), \v) - (p,\div \v) &= (\mathbf{f},\v) &\, \forall\, \v \in \Hozv,\\
	    -(\div \u, q) &= 0 &\, \forall\, q \in L^2(\Omega). \nonumber
    \end{alignat}
This corresponds to the standard Newton linearization of the Navier--Stokes equations with an augmented Lagrangian term, plus the linearization of the Lorentz force $\mathcal{D}$. We follow the approach of~\cite{Hong2015,Farrell2018,Farrell2020} to solve this system. The first idea is to use the augmented Lagrangian term $-\gamma \grad \div \u$ to approximate the inner Schur complement of the hydrodynamic block by choosing a large $\gamma$, e.g.~$\gamma \approx 10^4$. One can show \cite[Theorem 3.2]{Bacuta2006} that the inner Schur complement of the augmented system $\tilde{\mathcal{S}}_\text{NS}$ satisfies
\begin{equation}\label{eq:SchurCompNS}
	\tilde{\mathcal{S}}_\text{NS}^{-1} = \mathcal{S}_\text{NS}^{-1}-\gamma \MM_p^{-1}
\end{equation}
where $\mathcal{S}_\text{NS}^{-1}$ denotes the Schur complement of the system without the augmented Lagrangian term and $\MM_p$ denotes the pressure mass matrix. Therefore, for large $\gamma$ the pressure mass matrix scaled by $-1/(1/\Re + \gamma)$ is a good approximation for $\tilde{\mathcal{S}}_\text{NS}$. As the discretization considered in this work uses discontinuous pressures, the pressure mass matrix is block-diagonal and hence directly invertible. In the transient case $\mathcal{S}_\text{NS}^{-1}$ can be further approximated by the inverse of the stationary Schur complement plus an extra term $- \Delta t L_p^{-1}$  \cite{Heister2013}, where $L_p$ corresponds to the Poisson problem for $p$ with Neumann boundary conditions. In our numerical examples this extra term makes little difference as we only consider timesteps $\frac{1}{\Delta t} \ll \gamma$, and we therefore neglect it.

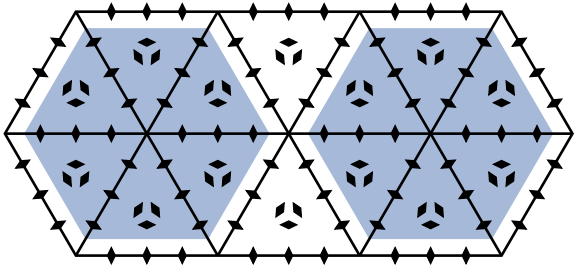
\begin{figure}[htbp]
	\centering
	\begin{tikzpicture}
  \foreach \i in {0, 1, 2, 3, 4} {
    \coordinate (0-\i) at ($(\i*1.5, 0)$);
    \coordinate (1-\i) at ($(0-\i) + (60:1.5)$);
    \coordinate (2-\i) at ($(1-\i) + (120:1.5)$);
    \coordinate (3-\i) at ($(\i*1.5, 2.6-2*1.3/3)$);
    \coordinate (4-\i) at ($(\i*1.5, 1.3-1*1.3/3)$);
    \coordinate (5-\i) at ($(0.75+ \i*1.5, 2.6-1*1.3/3)$);
    \coordinate (6-\i) at ($(0.75+ \i*1.5, 1.3-2*1.3/3)$);
  }
  
  \foreach \i in {1, 2, 3, 4}  {
    \node[transform shape,
    rotate=90,
    shape=diamond,
    fill=black,
    black,
    inner sep=0pt,
    minimum height=5pt,
    minimum width=2.5pt, draw] at ([shift={(90:3pt)}]4-\i) {};

    \node[transform shape,
    rotate=30,
    shape=diamond,
    fill=black,
    black,
    inner sep=0pt,
    minimum height=5pt,
    minimum width=2.5pt, draw] at ([shift={(30:-3pt)}]4-\i) {};

    \node[transform shape,
    rotate=-30,
    shape=diamond,
    fill=black,
    black,
    inner sep=0pt,
    minimum height=5pt,
    minimum width=2.5pt, draw] at ([shift={(-30:3pt)}]4-\i) {};
  }

  \foreach \i in {1, 2, 3, 4}  {
    \node[transform shape,
    rotate=90,
    shape=diamond,
    fill=black,
    black,
    inner sep=0pt,
    minimum height=5pt,
    minimum width=2.5pt, draw] at ([shift={(270:3pt)}]3-\i) {};

    \node[transform shape,
    rotate=210,
    shape=diamond,
    fill=black,
    black,
    inner sep=0pt,
    minimum height=5pt,
    minimum width=2.5pt, draw] at ([shift={(210:-3pt)}]3-\i) {};

    \node[transform shape,
    rotate=150,
    shape=diamond,
    fill=black,
    black,
    inner sep=0pt,
    minimum height=5pt,
    minimum width=2.5pt, draw] at ([shift={(150:3pt)}]3-\i) {};
  }
  
  \foreach \i in {1, 2, 3}  {
    \node[transform shape,
    rotate=90,
    shape=diamond,
    fill=black,
    black,
    inner sep=0pt,
    minimum height=5pt,
    minimum width=2.5pt, draw] at ([shift={(90:3pt)}]5-\i) {};

    \node[transform shape,
    rotate=30,
    shape=diamond,
    fill=black,
    black,
    inner sep=0pt,
    minimum height=5pt,
    minimum width=2.5pt, draw] at ([shift={(30:-3pt)}]5-\i) {};

    \node[transform shape,
    rotate=-30,
    shape=diamond,
    fill=black,
    black,
    inner sep=0pt,
    minimum height=5pt,
    minimum width=2.5pt, draw] at ([shift={(-30:3pt)}]5-\i) {};
  }

  \foreach \i in {1, 2, 3}  {
    \node[transform shape,
    rotate=90,
    shape=diamond,
    fill=black,
    black,
    inner sep=0pt,
    minimum height=5pt,
    minimum width=2.5pt, draw] at ([shift={(270:3pt)}]6-\i) {};

    \node[transform shape,
    rotate=210,
    shape=diamond,
    fill=black,
    black,
    inner sep=0pt,
    minimum height=5pt,
    minimum width=2.5pt, draw] at ([shift={(210:-3pt)}]6-\i) {};

    \node[transform shape,
    rotate=150,
    shape=diamond,
    fill=black,
    black,
    inner sep=0pt,
    minimum height=5pt,
    minimum width=2.5pt, draw] at ([shift={(150:3pt)}]6-\i) {};
  }

  \draw[line join=miter, thick]
  (0-1) -- (0-2) -- (0-3) -- (0-4)
  -- (1-4)
  -- (2-4)
  -- (2-3) -- (2-2) -- (2-1)
  -- (1-0)
  -- cycle;
  \draw[line join=miter, thick] (1-0) -- (1-4);
  \foreach \i/\j in {1/2, 2/3, 3/4} {
    \foreach \k in {0, 2} {
      \foreach \n in {0.25, 0.5, 0.75} {
        \node[transform shape,
        shape=diamond,
        fill=black,
        black,
        inner sep=0pt,
        minimum height=5pt,
        minimum width=2.5pt, draw] at ($(\k-\i)!\n!(\k-\j)$) {};
      }
    }
    \foreach \n in {0.25, 0.5, 0.75} {
      \node[transform shape,
      rotate=60,
      shape=diamond,
      fill=black,
      black,
      inner sep=0pt,
      minimum height=5pt,
      minimum width=2.5pt, draw] at ($(0-\i)!\n!(1-\i)$) {};

      \node[transform shape,
      rotate=-60,
      shape=diamond,
      fill=black,
      black,
      inner sep=0pt,
      minimum height=5pt,
      minimum width=2.5pt, draw] at ($(1-\i)!\n!(2-\i)$) {};

      \node[transform shape,
      rotate=-60,
      shape=diamond,
      fill=black,
      black,
      inner sep=0pt,
      minimum height=5pt,
      minimum width=2.5pt, draw] at ($(0-\j)!\n!(1-\i)$) {};

      \node[transform shape,
      rotate=60,
      shape=diamond,
      fill=black,
      black,
      inner sep=0pt,
      minimum height=5pt,
      minimum width=2.5pt, draw] at ($(1-\i)!\n!(2-\j)$) {};
    }
    \draw[line join=miter, thick] (0-\i) -- (1-\i);
    \draw[line join=miter, thick] (1-\i) -- (2-\i);
    \draw[line join=miter, thick] (0-\j) -- (1-\i);
    \draw[line join=miter, thick] (1-\i) -- (2-\j);
  }

  \foreach \n in {0.25, 0.5, 0.75} {
      \node[transform shape,
      rotate=-60,
      shape=diamond,
      fill=black,
      black,
      inner sep=0pt,
      minimum height=5pt,
      minimum width=2.5pt, draw] at ($(0-1)!\n!(1-0)$) {};

      \node[transform shape,
      rotate=60,
      shape=diamond,
      fill=black,
      black,
      inner sep=0pt,
      minimum height=5pt,
      minimum width=2.5pt, draw] at ($(1-0)!\n!(2-1)$) {};

      \node[transform shape,
      rotate=60,
      shape=diamond,
      fill=black,
      black,
      inner sep=0pt,
      minimum height=5pt,
      minimum width=2.5pt, draw] at ($(0-4)!\n!(1-4)$) {};

      \node[transform shape,
      rotate=-60,
      shape=diamond,
      fill=black,
      black,
      inner sep=0pt,
      minimum height=5pt,
      minimum width=2.5pt, draw] at ($(1-4)!\n!(2-4)$) {};

      \foreach \i/\j in {0/1, 1/2, 2/3, 3/4} {
        \node[transform shape,
        shape=diamond,
        fill=black,
        black,
        inner sep=0pt,
        minimum height=5pt,
        minimum width=2.5pt, draw] at ($(1-\i)!\n!(1-\j)$) {};
        }
  }

  \begin{scope}[on background layer]
    \draw[blue!50, fill=blue!50] 
    ([shift={(0:6pt)}]1-0) --
    ([shift={(60:6pt)}]0-1) --
    ([shift={(120:6pt)}]0-2) --
    ([shift={(180:6pt)}]1-2) --
    ([shift={(240:6pt)}]2-2) --
    ([shift={(300:6pt)}]2-1) --
    cycle;
    \draw[blue!50, fill=blue!50]
    ([shift={(0:6pt)}]1-2) --
    ([shift={(60:6pt)}]0-3) --
    ([shift={(120:6pt)}]0-4) --
    ([shift={(180:6pt)}]1-4) --
    ([shift={(240:6pt)}]2-4) --
    ([shift={(300:6pt)}]2-3) --
    cycle;
  \end{scope}
\node[fit=(current bounding box),inner sep=2mm]{};
\end{tikzpicture}
	\caption{Star patch for $\mathbb{BDM}_2$-elements.}
	\label{fig:star}
\end{figure}

Since the augmented Lagrangian term has a large kernel that consists of all solenoidal vector fields, a robust multigrid scheme as described in Section \ref{sec:robustssc} must be used to solve the augmented momentum block.
For the $\hdiv\times L^2$-conforming discretization the \emph{star iteration} \cite[section 4]{Farrell2018} can be used as a robust relaxation method. The subspace decomposition is defined as
\begin{equation}\label{eq:stardecomp}
	\mathbf{V}_i = \{\v \in \mathbf{V}_h: \mathrm{supp}(\v) \subset K_i \}
\end{equation}
where $K_i$ is the patch of elements sharing the vertex $i$ in the mesh. Example patches are shown in Figure~\ref{fig:star}. Since we use a structure-preserving discretization, the properties of the de Rham complexes \eqref{eq:deRhamdiscret3d} and \eqref{eq:deRhamdiscret2d} imply that \eqref{eq:stardecomp} fulfils the kernel decomposition property \eqref{eq:kerneldecomp}. This property was also used in \cite{Arnold2000} to construct a robust smoother for the $\hdiv$ and $\hcurl$ Riesz maps and in \cite{Hong2015} for the Stokes equations.
In this case we may employ the standard prolongation operator induced by the finite element discretization, because the uniformly-refined mesh hierarchy we consider is nested.

The velocity block further includes terms given by the convection-diffusion term $(\u\cdot \grad) \u$, the linearization of the Lorentz force $S\, \B^n \times (\u \times \B^n)$ and the stabilization term \eqref{eq:stabBurman}.
Numerical experiments in \cite{Farrell2020} and in the next Section~\ref{sec:numericalresults} show that these terms only degrade the performance of the preconditioner at high Reynolds and coupling numbers. As we have mentioned before, these somewhat surprising numerical observations are not backed up by theory since these terms do not fit in the framework of Section \ref{sec:robustssc}, and applying geometric multigrid methods to problems with strong advection typically requires special care.  The kernel of the stabilization $\mathcal{ST}(\u, \v)$ consists of all $C^1$ vector fields. Therefore, the stabilization term slightly degrades the performance of the solver, but the impact is not very significant as the factor $\mu h_{\partial K}^2$ is small.

\subsection{Solver for the electromagnetic block}\label{sec:solverformagneticblock}
The weak formulation of the electromagnetic block is given by
\begin{equation}
	\begin{aligned}
		(\E, \F) - \frac{1}{\Rem}(\B, \vcurl \F) + \delta\, (\u^n \times \B, \F) &= 0 &&\forall\, \F \in \hzcurl,\\
		\frac{\eta}{\Delta t} (\B, \C) + (\vcurl \E, \C) + \frac{1}{\Rem}(\div \B, \div \C) &= (\f, \C) &&\forall\, \C \in \hzdiv.
	\end{aligned}
\end{equation}
Recall that $\eta,\delta \in \{0,1\}$ distinguish between the stationary ($\eta=0$) and transient ($\eta = 1$) cases and the Picard ($\delta=0$) and Newton ($\delta = 1$) linearizations.
Eliminating $\E$, this corresponds to a mixed formulation of
\begin{equation}
	\begin{aligned}
		\frac{\eta}{\Delta t} \B +  \frac{1}{\Rem}\left(\vcurl \vcurl \B - \grad \div \B\right)   + \delta\, \vcurl(\u^n\times\B) &= \f \text{ in } \Omega,\label{eq:LaplaceGtilde} \\
		\B \cdot \mathbf{n} &= 0 \text{ on } \partial \Omega, \\
		\frac{1}{\Rem} \vcurl \B - \delta\, \u^n\times \B &= \mathbf{0} \text{ on } \partial \Omega.
	\end{aligned}
\end{equation}
 For the Picard linearization, this problem simplifies to the mixed formulation for the standard vector Laplace problem \cite{Arnold2006} with boundary conditions $\B\cdot \n = \vcurl \B =\mathbf{0}$ on $\del\Omega$. Chen et al.\ \cite{Chen2018} propose a Schur complement solver and Arnold et al.\ \cite[\S 10]{Arnold2006} propose a norm-equivalent block diagonal preconditioner for the mixed formulation.
We also found that the star multigrid solver applied monolithically to the electromagnetic block \eqref{eq:LaplaceGtilde} results in an efficient solver and employ this solver in our numerical examples.
All of the solvers described show $\Rem$-robust performance.

In contrast, the presence of the additional term $\vcurl(\u^n \times \B)$ in the Newton linearization, which has a non-trivial kernel, makes the problem almost singular for high $\Rem$ in the stationary case and hence requires a special multigrid method.
Unfortunately the troublesome term $\vcurl(\u^n \times \B)$ is not symmetric and thus does not fit the available analytical framework of Sch\"oberl. Our considerations on this point are therefore necessarily heuristic. Some insight may be gained by employing the vector identity
\begin{equation}
	\vcurl (\mathbf{A} \times \mathbf{B}) =  (\mathbf{B} \cdot \nabla) \mathbf{A} -  (\mathbf{A} \cdot \nabla) \mathbf{B} +
	\mathbf{A}(\nabla \cdot \mathbf{B}) - \mathbf{B}(\nabla \cdot \mathbf{A})
\end{equation}
to rewrite \eqref{eq:LaplaceGtilde} to
\begin{equation}
	\frac{\eta}{\Delta t} \B -\frac{1}{\Rem}\Delta \B - (\B \cdot \grad) \u^n + (\u^n \cdot \nabla) \B  - \u^n(\nabla \cdot \mathbf{B}) - \mathbf{B}(\nabla \cdot \u^n).
\end{equation}
The last term $ - \mathbf{B}(\nabla \cdot  \u^n)$ vanishes since we exactly enforce $\nabla \cdot \u^n=0$ in each step.
The terms  $-(\B \cdot \grad) \u^n + (\u^n \cdot \nabla) \B $ are reminiscent of the Newton linearization of the advection term $(\u \cdot \grad) \u$ of the Navier--Stokes equation, for which it has been demonstrated that a star multigrid method is effective \cite{Farrell2018}. Numerical experiments with these approaches applied monolithically do indeed yield a robust solver for the stationary and transient cases in two dimensions, and in the transient case in three dimensions for sufficiently small $\Delta t$. We observe in our numerical tests that in three dimensions the solver breaks down for $\Rem \approx 700$ for a stationary lid-driven cavity problem.

\section{Numerical results} \label{sec:numericalresults}

In this section, we present numerical results for the Picard and Newton linearization described in the previous sections. We investigate three test problems: the stationary and transient version of a lid-driven cavity problem and a transient island-coalescence problem. The numerical results were produced on ARCHER2, the UK national supercomputer, which consists of 5,860 compute nodes each built of two AMD Zen2 7742 processors with 64 2.25GHz cores  and 256GB of memory.

\subsection{Algorithm details}\label{sec:algorithmdetails}
The algorithm is implemented in Firedrake \cite{rathgeber2016} and uses the solver packages PETSc \cite{balay2019} and PCPATCH \cite{farrell2019pcpatch}. The latter includes the implementation of the multigrid relaxation method described in Section \ref{sec:derivationofblockpreconditioner}. It is well-known that the convergence of the nonlinear scheme depends heavily on the initial guess and might fail to converge for high Reynolds numbers with poor initial guesses. To circumvent this problem we perform continuation in the Reynolds numbers and coupling number, for the stationary problems. In the presented tables we always apply continuation to the variable in the column first and use each solution as the starting point for the continuation over the rows. We use the steps $1,100,1000,5000,10000$ for $S$ and $1,500,1000,3000,5000,7000,9000,10000$ for $\Re$ and $\Rem$.
The reported nonlinear and linear iteration numbers correspond to the final solve in the continuation; however, the extra cost for the continuation should be kept in mind for stationary problems. For time-dependent problems, continuation is not necessary.

We use flexible GMRES \cite{saad1993} as the outermost Krylov solver since we apply GMRES in the multigrid relaxation. Moreover, we apply a block upper triangular preconditioner \cite{Benzi2005}
\begin{equation}
	\mathcal{P}=
	\begin{pmatrix}
		\II &-\tilde{\MM}^{-1}\mathcal{K} \\
		\mathbf{0} &\II
	\end{pmatrix}
	\begin{pmatrix}
		\tilde{\MM}^{-1} &\mathbf{0} \\
		\mathbf{0} &\tilde{\mathcal{S}}^{-1}
	\end{pmatrix}
\end{equation}
to \eqref{eq:matrix_upBE}, where we denoted \eqref{eq:matrix_upBE} here as $\begin{pmatrix} \MM & \mathcal{K} \\ \mathcal{L} & \mathcal{N} \end{pmatrix}$. We also investigated a full block-LDU preconditioner without notable improvements in terms of iteration counts, which fits with the recent theoretical results in \cite{Southworth2020}.

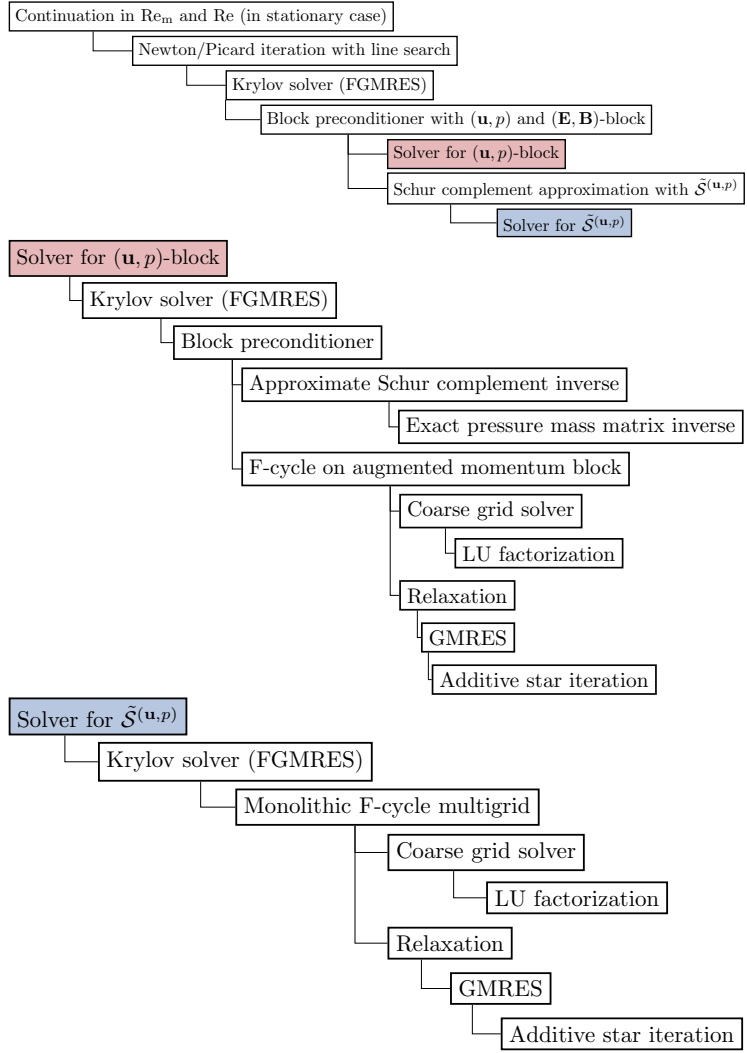
\begin{figure}[htb!]
	\centering
	\resizebox{0.75\textwidth}{!}{\begin{tikzpicture}[%
  every node/.style={draw=black, thick, anchor=west},
  grow via three points={one child at (-1.4,-0.7) and
  two children at (-1.4,-0.7) and (-1.4,-1.4)},
  edge from parent path={(\tikzparentnode.188) |- (\tikzchildnode.west)}]
  \node {Continuation in $\Rem$ and $\Re$ (in stationary case)}
  child { node{Newton/Picard iteration with line search}
  child { node{Krylov solver (FGMRES)}
    child { node {Block preconditioner with $(\u, p)$ and $(\E,\B)$-block}
      child { node[fill=red!40] {Solver for $(\u, p)$-block}
      }
      child { node {Schur complement approximation with $\tilde{\mathcal{S}}^{(\u, p)}$}
      child { node[fill=blue!40] {Solver for $\tilde{\mathcal{S}}^{(\u, p)}$}
      }
      }
     }
     }
    };
    
\end{tikzpicture}}
	\resizebox{0.75\textwidth}{!}{\begin{tikzpicture}[%
  every node/.style={draw=black, thick, anchor=west},
  grow via three points={one child at (-0.6,-0.7) and
  two children at (-0.6,-0.7) and (-0.6,-1.4)},
  edge from parent path={(\tikzparentnode.201) |- (\tikzchildnode.west)}]
   \node[fill=red!40]{Solver for $(\u, p)$-block}
      child { node {Krylov solver (FGMRES)}
        child { node {Block preconditioner}
          child { node {Approximate Schur complement inverse}
              child{ node {Exact pressure mass matrix inverse}}
          }
          child [missing] {}
          child { node {F-cycle on augmented momentum block}
              child { node {Coarse grid solver}
                child { node {LU factorization}}
              }
              child [missing] {}
              child { node {Relaxation}
                child { node {GMRES}
                  child { node {Additive star iteration}}
                }
              }
          }
        }
      };
\end{tikzpicture}}
	\resizebox{0.75\textwidth}{!}{\begin{tikzpicture}[%
  every node/.style={draw=black, thick, anchor=west},
  grow via three points={one child at (0.0,-0.7) and
  two children at (0.0,-0.7) and (0.0,-1.4)},
  edge from parent path={(\tikzparentnode.210) |- (\tikzchildnode.west)}]
   \node[fill=blue!40] {Solver for $\tilde{\mathcal{S}}^{(\u, p)}$}
      child { node {Krylov solver (FGMRES)}
          child { node {Monolithic F-cycle multigrid}
              child { node {Coarse grid solver}
                child { node {LU factorization}}
              }
              child [missing] {}
              child { node {Relaxation}
                child { node {GMRES}
                  child { node {Additive star iteration}}
                }
              }
          }
      };
\end{tikzpicture}}
	\caption{Graphical outline of the solver.}
	\label{fig:solver_diagramm}
\end{figure}

Both the block matrix $\MM$ and the outer Schur complement approximation $\tilde{\mathcal{S}}^{(\u, p)}$ are inverted approximately with two iterations of preconditioned FGMRES (denoted $\tilde{\MM}^{-1}$ and $\tilde{S}^{-1}$ respectively). The former uses the block preconditioner for the hydrodynamic block described in Section \ref{sec:solverforschurcomp}, the latter the monolithic multigrid method described in Section \ref{sec:solverformagneticblock}. In the numerical results we focus on taking the outer Schur complement that eliminates the hydrodynamic block, except for one case in Section \ref{sec:stationaryliddrivencavityproblemin3d}. In both multigrid methods we use six preconditioned GMRES iterations as the smoother on each level and the direct solver MUMPS \cite{MUMPS} to solve the problem on the coarsest grid. Since this relaxation is quite expensive, convergence in a very small number of outer iterations is required for efficiency. See Figure \ref{fig:solver_diagramm} for a graphical representation of the solver.

We have chosen relative and absolute tolerances of $10^{-10}$ and $10^{-6}$ for the nonlinear solver and $10^{-7}$ and $10^{-7}$ for the outermost linear solver, measured in the Euclidean norm. We use the  $\hdiv\times L^2$-conforming elements $\mathbb{BDM}_2 \times \mathbb{DG}_{1}$ for $(\u_h,p_h)$.   Moreover, we apply $\mathbb{CG}_2 \times \mathbb{RT}_2$ elements for $(E_h,\B_h)$ in 2D and $\mathbb{NED}1_2 \times \mathbb{RT}_2$ elements for $(\E_h,\B_h)$ in 3D. When we consider a manufactured solution we always subtract $\int_\Omega p \,\mathrm{d}x$ from the pressure to fix the average of $p$ to be zero.

For time-dependent problems, we apply the second-order, L-stable BDF2 method with a fixed time-step. We compute the first time-step with Crank-Nicolson to provide the second starting value for BDF2. For the transient lid-driven cavity problem, we use a time-step of $\Delta t = 0.01$ and a final time of $T=0.1$. We did not choose a higher final time T because of budget limitations. However, we confirmed that the reported numbers are representative for higher T by computing the solution for a few parameters until $T=1$ without noticeable changes in the iteration counts. Moreover, we confirm the efficiency for more time-steps in the island coalescence problem where we iterate in the finest run until $T=15$ in 2400 timesteps.


\subsection{Interpolating boundary data}
\label{sec:interp-boundary-data}
The theory from the previous sections has been formulated for homogeneous boundary conditions, but the generalisation is straightforward for non-homogeneous boundary conditions. However, there is a subtle technicality in the implementation if one wants to enforce the divergence constraint $\div \B_h=0$ pointwise. Strong boundary conditions are enforced in a finite element code by interpolating the given boundary data onto the corresponding finite element space. If the interpolation of the boundary values $\mathbf{g}$ were exact, identity \eqref{eq:InterpolationPreserveRT} would imply that $\div \B_h=0$ holds. However, the degrees of freedom for the interpolation are moments and are usually implemented by a quadrature rule whose quadrature degree is based on the polynomial degree of the finite element space. If $\mathbf{g}$ is a non-polynomial expression, this quadrature rule might not interpolate the boundary condition exactly and therefore one loses the property that $\div \mathbf{g}_h=0$ on $\del \Omega$ holds exactly.

\begin{figure}[htb!]
	\centering
	\includegraphics[width=0.5\textwidth]{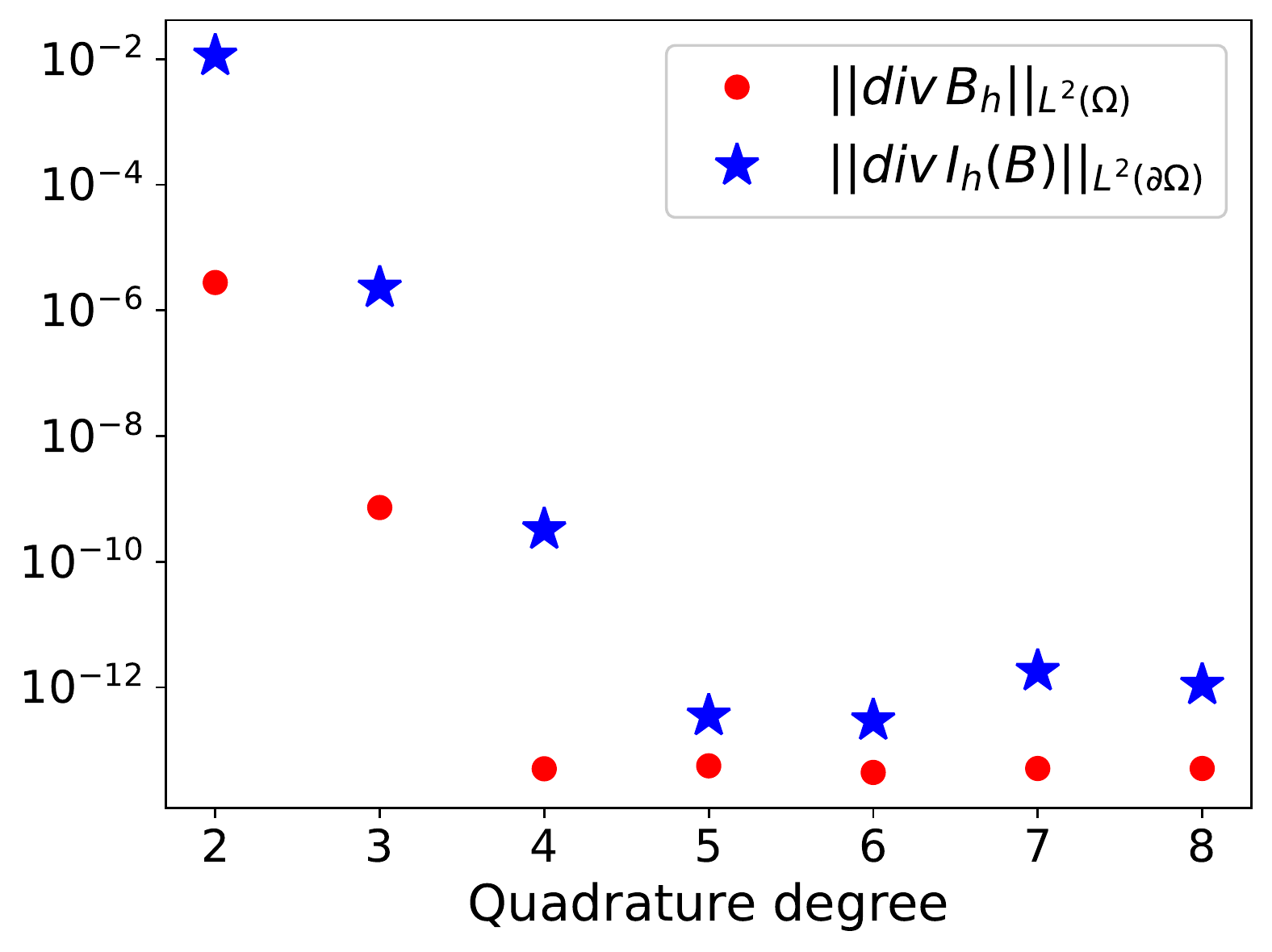}
	\caption{$L^2$-norm of the divergence of the solution $\B_h$ and the interpolant of the boundary condition for different quadrature degrees in the evaluation of the degrees of freedom for the Raviart--Thomas space.}
	\label{tab:divBQuadDeg}
\end{figure}

To circumvent this problem we use high-order quadrature rules for the evaluation of the degrees of freedom to ensure that the interpolation is exact up to machine precision.
In Figure \ref{tab:divBQuadDeg} we have illustrated the effect of the quadrature degree on the enforcement of the divergence constraint. We have used the method of manufactured solutions for a smooth problem to compute $\|\div \B_h\|_0$ for different quadrature degrees. Moreover, we have plotted the $L^2$-norm over $\partial \Omega$ of the interpolation of the divergence-free function $\B$ into the $\mathbb{RT}_2$ space. One can clearly observe that a quadrature degree of 2 for $\mathbb{RT}_2$ elements is not sufficient to enforce $\div \B_h=0$ up to machine precision. A higher quadrature degree preserves the divergence of the boundary data more accurately and leads to the point-wise enforcement of $\div \B_h=0$.

\subsection{Two-dimensional results}

\subsubsection{Stationary lid-driven cavity in two dimensions}

First, we consider a lid-driven cavity problem posed over $\Omega=(-1/2,1/2)^2$ for a background magnetic field $\B_0=(0,1)^\top$ which determines the boundary conditions $\B\cdot \n = \B_0\cdot \n$ on $\del \Omega$, and set $\f = \mathbf{0}$ \cite{Ma2016}. We impose the boundary condition $\u=(1,0)^\top$ at the boundary $y=0.5$ and homogeneous boundary conditions elsewhere. The problem models the flow of a conducting fluid driven by the movement of the lid at the top of the cavity. The magnetic field imposed orthogonal to the lid creates a Lorentz force that perturbs the flow of the fluid. 

For the multigrid hierarchy we use a coarse mesh of $16 \times 16$ cells and six levels of refinement in 2D resulting in a $1024\times1024$ grid with 73.4 million degrees of freedom (DoFs). For both linearizations we observe fairly constant Krylov iteration counts for $\Re$ and $S$ in the range of 1 to 10,000 in Table \ref{tab:ldc2dSRE}. In terms of the nonlinear convergence, the Picard linearization takes sometimes slightly more iterations than the Newton linearizations, with slightly better linear iteration numbers.

\begin{table}[htbp!]
	\centering
	\begin{tabular}{r|ccc|ccc}
		\toprule
		& \multicolumn{3}{c|}{Picard} &  \multicolumn{3}{c}{Newton} \\
		\midrule
		$S\backslash\Re$ & 1        & 1,000    & 10,000     & 1        & 1,000    & 10,000   \\
		\midrule
		1           & ( 3) 5.3 & ( 4) 3.5 & ( 3) 4.3 & ( 2) 6.5 & ( 4) 3.5 & ( 3) 4.3   \\
		1,000     & ( 4) 3.5 & ( 3) 4.7 & ( 2) 8.5 & ( 2) 5.5 & ( 3) 4.7 & ( 2) 6.5     \\
		10,000    & ( 3) 5.0 & ( 3) 4.3 & ( 2) 7.0 & ( 2) 6.5 & ( 2) 6.0 & ( 2) 7.0     \\
		\bottomrule
	\end{tabular}
	\caption{(Nonlinear iterations) Average outer Krylov iterations per nonlinear step for the stationary lid-driven cavity problem in 2D.\label{tab:ldc2dSRE}}
\end{table}

\begin{figure}[htbp!]
	\centering
	\begin{tabular}{ccc}
		\includegraphics[width=3.9cm]{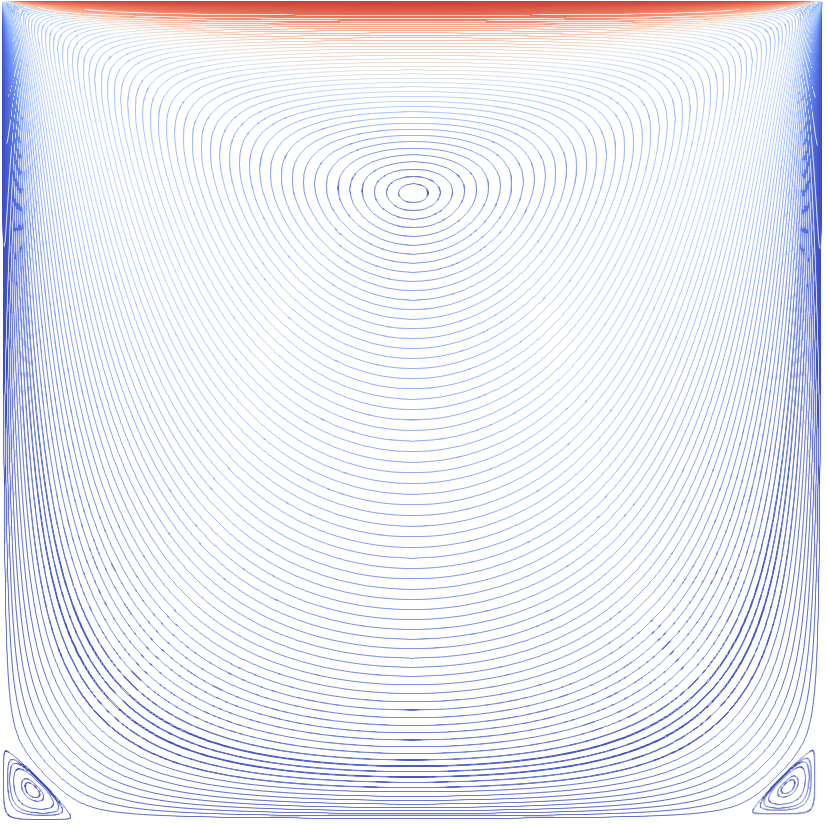} &
		\includegraphics[width=3.9cm]{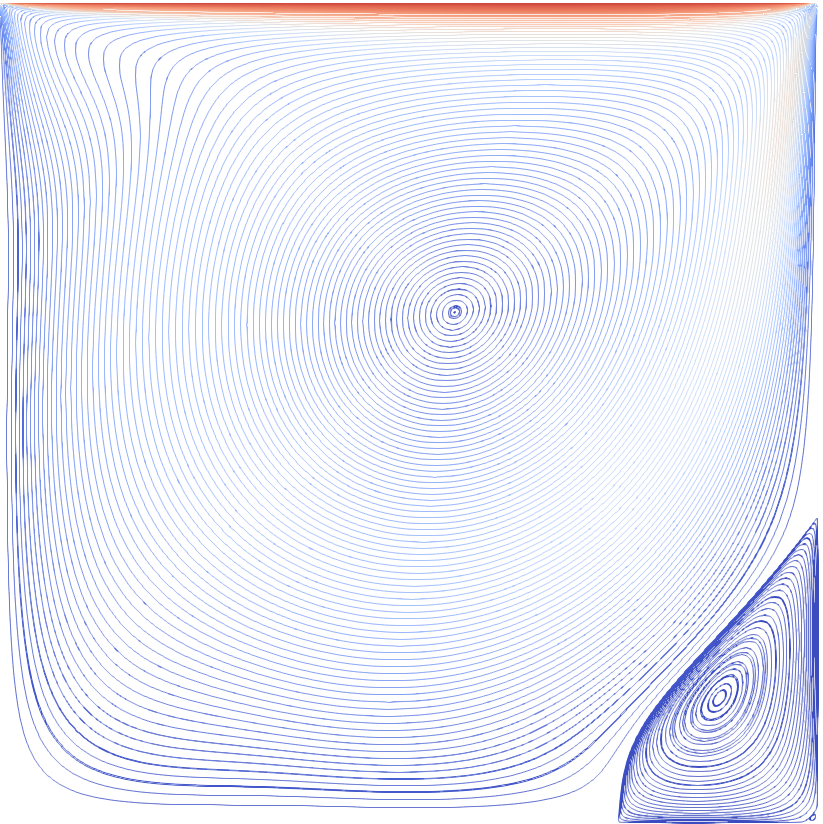} &
		\includegraphics[width=3.9cm]{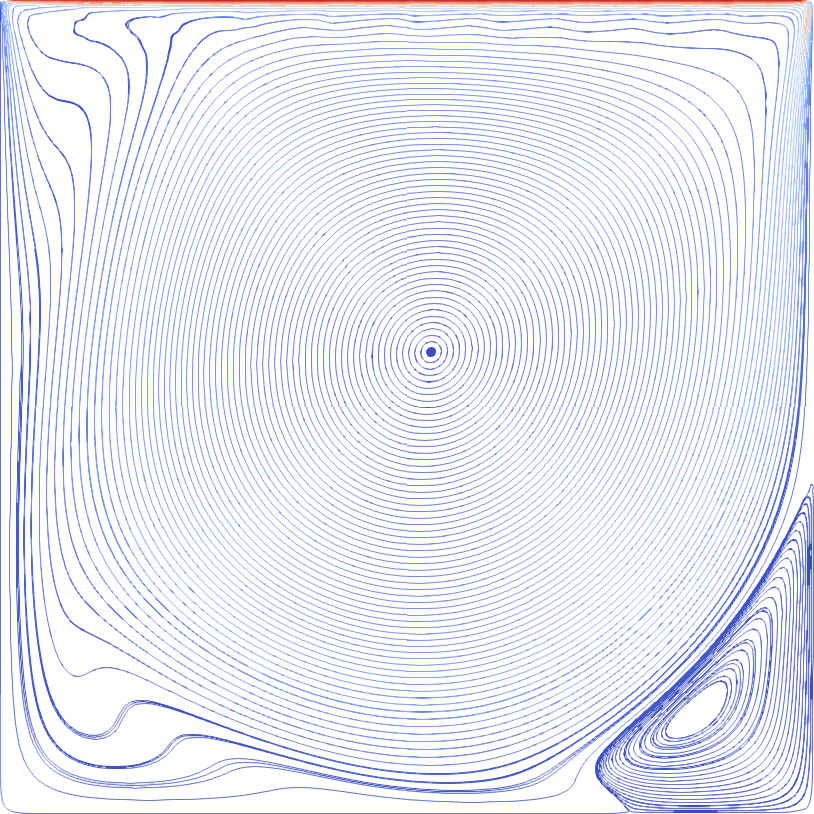} \\
		\includegraphics[width=3.9cm]{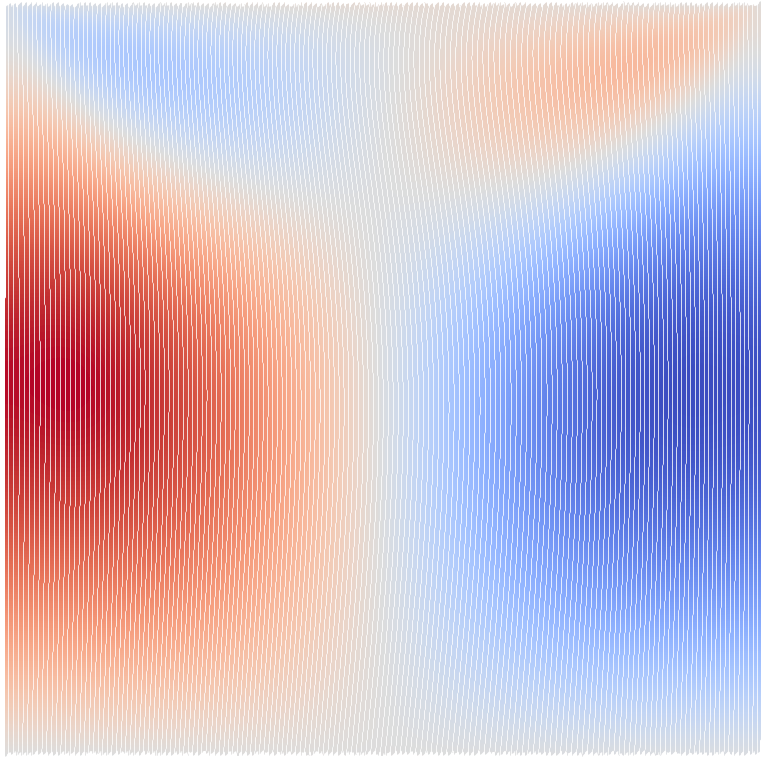} &
		\includegraphics[width=3.82cm]{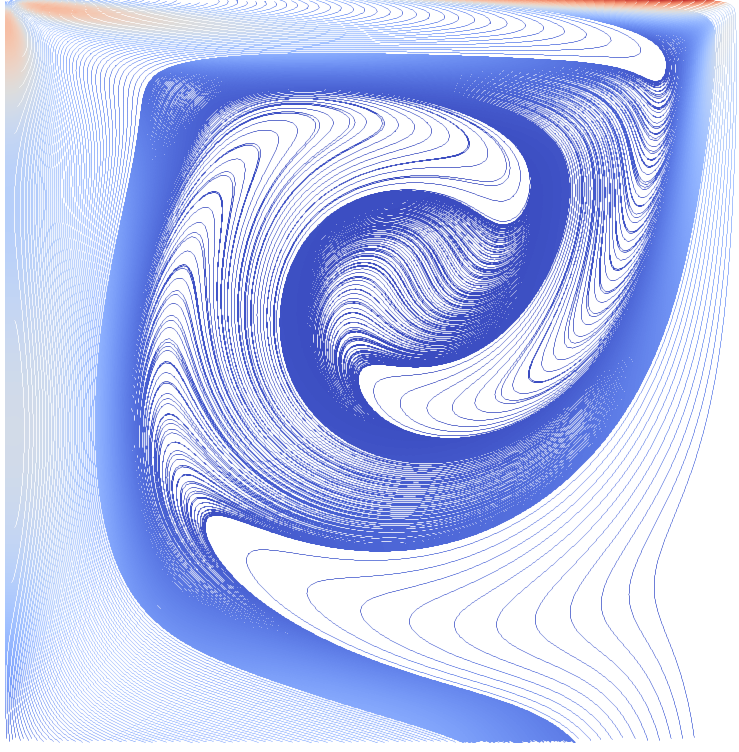} &
		\includegraphics[width=3.82cm]{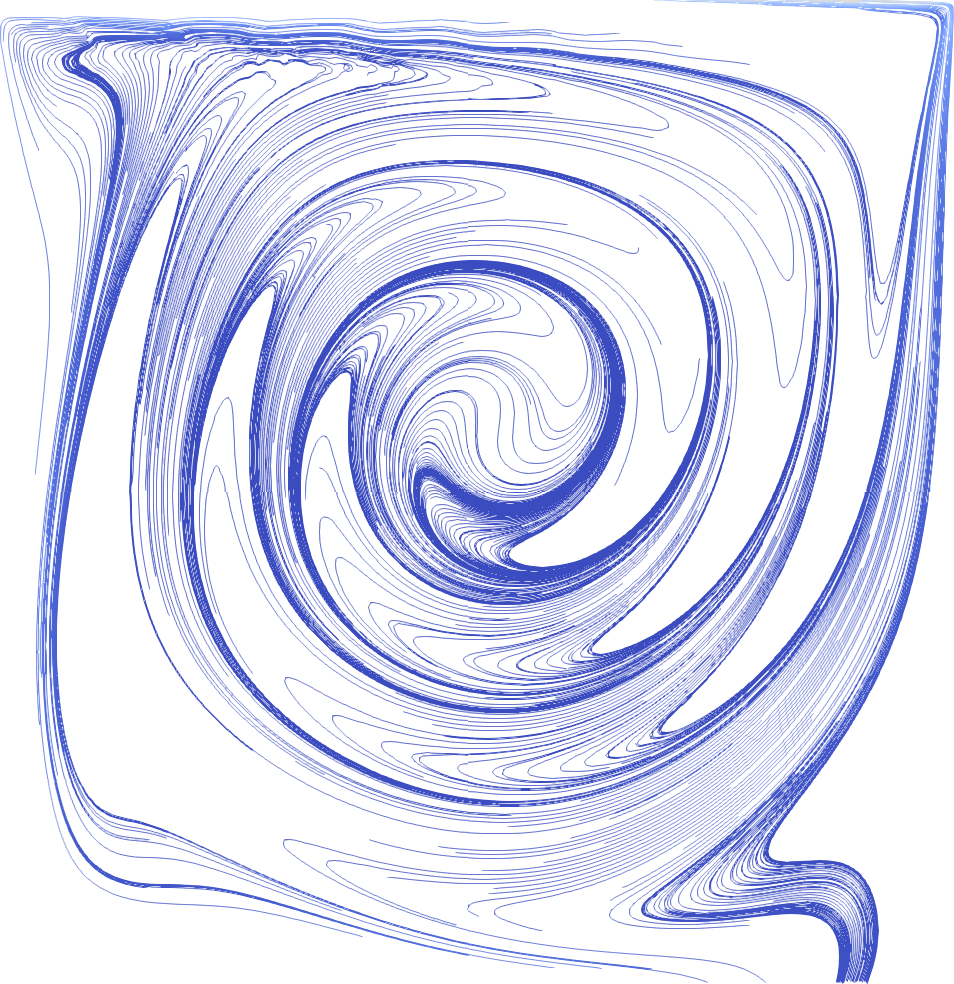} \\
		$\Re=\Rem=1$ & $\Re=\Rem=500$ & $\Re=\Rem=5,$000\\
	\end{tabular}
	\caption{Streamlines for the two-dimensional stationary lid-driven cavity problem for $\u$ (upper row) and $\B$ (lower row).\label{fig:Streamlinesldc}}
\end{figure}

As mentioned earlier, our scheme does not include a stabilization for high magnetic Reynolds numbers. However, we have verified that our solutions do not exhibit oscillations in this regime. A plot of the streamlines for different $\Re$ and $\Rem$ can be found in Figure \ref{fig:Streamlinesldc}. One can clearly observe the phenomenon that for high magnetic Reynolds numbers the magnetic field lines are advected with the fluid flow. Iteration counts for varying $\Re$ and $\Rem$ are displayed in Table \ref{tab:ldc2dREREMHdiv}. 

For the Picard linearization we observe that the nonlinear scheme already fails to converge for a magnetic Reynolds number of 100. The poor nonlinear convergence of the Picard iteration for high $\Rem$ even with continuation was previously observed for other formulations \cite{PhillipsPHD,Phillips2014}.


For the Newton linearization the linear iterations increase slightly since the approximation of the Schur complement $\mathcal{S}^{(\u, p)}_\text{N}$ by $\tilde{\mathcal{S}}^{(\u, p)}$ becomes less accurate for high Reynolds numbers. On the other hand, the number of nonlinear iterations remains fairly constant which seems to indicates that the linear solver for the $(\E_h, \B_h)$ block described in Section \ref{sec:solverformagneticblock} works very well for high $\Rem$ in two dimensions.


\begin{table}[htbp!]
  \centering
    \begin{tabular}{r|ccc|ccc}
      \toprule
      & \multicolumn{3}{c|}{Picard} & \multicolumn{3}{c}{Newton} \\
      \midrule
      $\Rem\backslash\Re$  &1 &     1,000 &    10,000 &1 &     1,000 &    10,000 \\
      \midrule
      1              &( 3) 5.3&(4) 3.5 & (3) 4.3&( 2) 6.0 & ( 3) 4.3 & ( 3) 4.3   \\
1,000      &-&-&-& ( 2) 4.5 & ( 3) 3.0 & ( 3) 3.0  \\
10,000    &-&-&-& ( 2) 4.5 & ( 4) 5.5 & ( 3) 5.7      \\
      \bottomrule
    \end{tabular}
  \caption{Iteration counts for the stationary lid-driven cavity problem in 2D with $\hdiv\times L^2$-discretization for different $\Rem$ and $\Re$.\label{tab:ldc2dREREMHdiv}}
    \begin{tabular}{r|ccc|ccc}
      \toprule
      & \multicolumn{3}{c|}{Picard} & \multicolumn{3}{c}{Newton} \\
      \midrule
      $\Rem\backslash\Re$ & 1        & 1,000    & 10,000  & 1        & 1,000    & 10,000   \\
      \midrule
      1                   & ( 2) 4.0 & ( 2) 2.5 & ( 2) 8.5 & ( 2) 4.0 & ( 2) 2.5 & ( 3) 9.7\\
      1,000               & -        & -        & -        & ( 5) 1.8 & ( 3) 3.0 & ( 2) 4.0 \\
      10,000              & -        & -        & -        & ( 8) 5.2 & ( 4) 6.2 & ( 2) 5.5 \\
      \bottomrule
    \end{tabular}
  \caption{Iteration counts for the stationary lid-driven cavity problem in 2D with Scott--Vogelius elements for different $\Rem$ and $\Re$.\label{tab:ldc2dREREMSV}}
\end{table}

\subsubsection{Scott--Vogelius discretization for $(\u_h, p_h)$}
Thus far we have considered a $\hdiv\times L^2$ discretization for the hydrodynamic variables. For comparison, in this subsection we include results for Scott--Vogelius elements \cite{scott_conforming_1985}, i.e.\ $(\mathbb{CG}_k)^d\times \mathbb{DG}_{k-1}$ elements. A fluid-Reynolds-robust preconditioner for this element pair was recently developed in \cite{Farrell2020}. While this conforming discretization does not require stabilization terms to weakly enforce continuity, it is only stable on certain types of meshes. For this reason, the mesh hierarchy is barycentrically refined, and the specialized multigrid method of \cite{Farrell2020} exploits this structure. This barycentric refinement ensures stability for polynomial order $k=d$ \cite{zhang2004}.

The results are shown in Table \ref{tab:ldc2dREREMSV}. We observe that the Krylov iteration counts are in general similar for the Scott--Vogelius element, making this an attractive alternative for those wishing to employ conforming schemes. However, one must keep in mind that the work per Krylov iteration is substantially higher for this element, due to the use of larger patches in the so-called macrostar relaxation method.

To summarize the two-dimensional stationary results, both schemes considered provide $\Re$-robust solvers and also perform very well for large coupling numbers $S$. The Picard iteration is unsuitable for high magnetic Reynolds number because the nonlinear iteration fails to converge. The Newton scheme performs well for high $\Rem$ with good nonlinear convergence and a slight increase in the linear iteration numbers.

\subsubsection{Time-dependent lid-driven cavity problem in two dimensions}
We next consider the time-dependent lid-driven cavity problem. We choose the same boundary conditions and right-hand side as in the stationary case. The numerical results for varying $S$ and $\Re$ are shown in Table \ref{tab:ldc2dtimedep}. As in the stationary case, the Krylov iteration counts remain almost constant for the two linearizations. We notice that the Picard iteration fails to converge for high $S$ and $\Rem$ for the chosen $\Delta t = 0.01$. However, we tested that one can get the Picard iteration to converge in most cases by choosing a smaller $\Delta t$ in the first timesteps. We do not report these results here to keep the tables consistent.

Table \ref{tab:ldc2dtimedep} also shows iteration counts for varying $\Re$ and $\Rem$. The linear solver is robust for most parameter values, with iteration counts only increasing for $\Re=1$ and $\Rem=100,000$.

For completeness, we also study the case of high $\Rem$ and $S$ at the same time in Table \ref{tab:ldc2dtimedep}, which we expect to be the most challenging case. Again the slight increase of the Krylov iterations in the Newton iteration is due to inaccurate outer Schur complement approximation. However, the solvers perform very well, considering the difficulty of the problem.

\begin{table}[htbp!]
  \centering
    \begin{tabular}{r|ccc|ccc}
      \toprule
      & \multicolumn{3}{c|}{Picard} & \multicolumn{3}{c}{Newton} \\
      \midrule
      $S \backslash\Re$  &1 &     10,000 &    100,000 & 1 &     10,000 &    100,000 \\
      \midrule
      1 &  (2.0) 3.0 & (2.2) 3.6 & (3.1) 3.3 & (1.6) 3.6 & (2.2) 3.6 & (3.1) 3.3\\
      1,000 &   (3.0) 4.0 & (3.0) 3.3 & (2.7) 3.0 & (2.1) 4.7 & (2.2) 3.9 & (2.3) 3.3\\
      10,000 & - & - & - &  (2.2) 6.5 & (2.5) 5.0 & (2.3) 5.6\\
      \bottomrule
      \multicolumn{7}{c}{} \\
      \toprule
      & \multicolumn{3}{c|}{Picard} & \multicolumn{3}{c}{Newton} \\
      \midrule
      $\Rem\backslash\Re$  &1 &     10,000 &    100,000 &1 &     10,000 &    100,000 \\
      \midrule
      1 & (2.0) 3.0 & (2.2) 3.6 & (3.1) 3.3 & (1.6) 3.6 & (2.2) 3.6 & (3.1) 3.3  \\
      10,000 & - & - & - & (2.0) 3.1 & (2.3) 3.6 & (3.1) 3.3\\
      100,000 &  - & - & - &  (2.2)10.9 & (3.0) 3.2 & (3.3) 3.2\\
      \bottomrule
      \multicolumn{7}{c}{} \\
      \toprule
      & \multicolumn{3}{c|}{Picard} & \multicolumn{3}{c}{Newton} \\
      \midrule
      $\Rem\backslash S$  &1 &     1,000 &    10,000 & 1 &     1,000 &    10,000 \\
      \midrule
      1 &  (2.0) 3.0 & (3.0) 4.0 & - & (1.6) 3.6 & (2.1) 4.7 & (2.2)  6.5  \\
      1,000 &  (3.0) 2.5 & - & - & (2.0) 3.1 &  (2.2) 5.6 & (2.8)11.0 \\
      10,000 & - & - & - &  (2.0) 3.1 & (2.2) 6.3 & (3.1)11.8\\
      \bottomrule
    \end{tabular}
  \caption{Iteration counts for the transient lid-driven cavity problem in 2D.\label{tab:ldc2dtimedep}}
\end{table}

\subsubsection{Time-dependent island coalescence problem in two dimensions}
Next, we consider a two-dimensional island coalescence problem to demonstrate the effectiveness of our method for a physically relevant model that shows behaviour which is unique to MHD problems. Furthermore, we report results for a weak scalability test with up to 256 processors and 18.3M DoFs to examine the performance of our algorithm. \\
The island coalescence problem is used to model magnetic reconnection processes in large aspect ratio tokamaks. For a strong magnetic field in the toroidal direction, the flow can be described in a two-dimensional setting by considering a cross-section of the tokamak. We consider the same problem as in \cite[Section 4.2]{adler2020monolithic}. The domain $\Omega=(-1,1)^2$ results from the unfolding of an annulus in the cross-sectional direction where the left and right edges are mapped periodically. The equilibrium solution for $k=0.2$ is given by
\begin{align*}
 \u_{eq}&=\mathbf{0},\qquad p_{eq}(x,y) = \frac{1-k^2}{2}\left (1+\frac{1}{(\cosh(2\pi y) + k \cos(2\pi x))^2}\right), \\
\B_{eq}(x,y) &= \frac{1}{\cosh (2\pi y) + k \cos(2\pi x)} 
\begin{pmatrix}
	\sinh (2\pi y) \\ k \sin(2\pi x)
\end{pmatrix}, E_{eq}=\frac{1}{\Rem} \scurl \B_{eq} - u_{eq} \times B_{eq},
\end{align*}
which results in right-hand sides $\f=\mathbf{0}$ and $\mathbf{g}$ for \eqref{eq:MHD3d3} given by
\begin{equation}
	\mathbf{g} = \frac{-8 \pi^2 (k^2-1)}{\Rem(\cosh(2\pi y) + k\cos(2\pi x))^3}
	 \begin{pmatrix}
		\sinh(2\pi y) \\ k \sin(2\pi x)
	\end{pmatrix}.
\end{equation}
The initial condition for $\B_{eq}$ is given by perturbing it for $\varepsilon=0.01$ with 
\begin{equation}
	\Delta \B = \frac{\varepsilon}{\pi} \begin{pmatrix}
		-\cos(\pi x) \sin(\pi y/2) \\
		2\cos(\pi y/2) \sin(\pi x)
	\end{pmatrix}.
\end{equation}
The authors believe that the reported $\Delta \B$ in \cite{adler2020monolithic} includes a typographical error, as it is not divergence-free, and amended the second component appropriately. Therefore, the problem setup is not exactly identical and hence we might see slightly different solutions.
The reconnection rate can be computed as the difference between $\curl \B$ evaluated at $(0,0)$ at the current time and the initial time, divided by $\sqrt{\Rem}$. Note that in our formulation $\B \in \hdivz$ and therefore we apply the $\scurl$ weakly by solving a problem for $j_0 \in \hzcurl$ such that
\begin{equation}
	(j_0, k) = (\B, \vcurl k) \quad \forall\, k \in \hzcurl.
\end{equation}
In order to make the point evaluation of $j_0$ at (0,0) well-posed we project $j_0$ to the space $\mathbb{CG}1$ as in \cite{adler2020monolithic}. 

Figure 4 shows the reconnection rates for $\Re=\Rem=1,000$, $\Re=\Rem=5,000$ and $\Re=\Rem=10,000$ for three different spatial and temporal resolutions. We have fixed a coarse grid of $16\times16$ cells and compute results for three (1.1M DoFs), four (4.6M DoFs) and five (18.4M DoFs) levels of refinement. For the three levels of refinement, we chose $\Delta t =0.025$ and halved it with each refinement. We iterated until $T=15$ which results in 2400 timesteps for the finest resolution.

 One can observe a decreasing height of the peak for increasing Reynolds numbers and the so-called ``sloshing" \cite{Knoll2006} effect that results in further peaks after the main peak with higher Reynolds numbers. Convergence for our considered meshes can be observed for $\Re=\Rem=1,000$ and  $\Re=\Rem=5,000$ while a further refinement is needed for $\Re=\Rem=10,000$. Nevertheless, our finest grid results match the results of \cite[Fig. 4]{adler2020monolithic} where finer meshes of up to $2560\times2560$ cells and $\Delta t=0.0016$ have been considered. For example, our finest result for $\Re=\Rem=10,000$ clearly reproduces the second peak in the reconnection rate.
 
 \begin{figure}[htbp!]
 	\centering
 	\begin{tabular}{ccc}
 		\includegraphics[width=3.9cm]{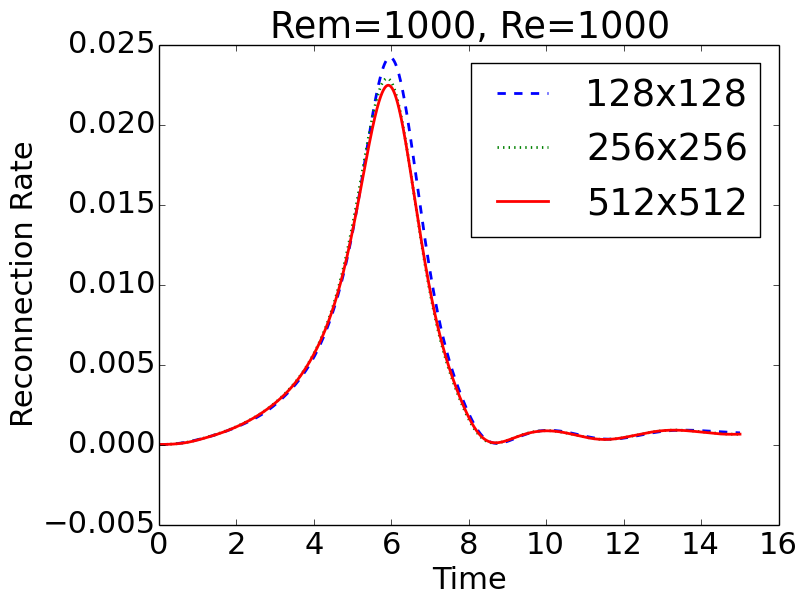} &
 		\includegraphics[width=3.9cm]{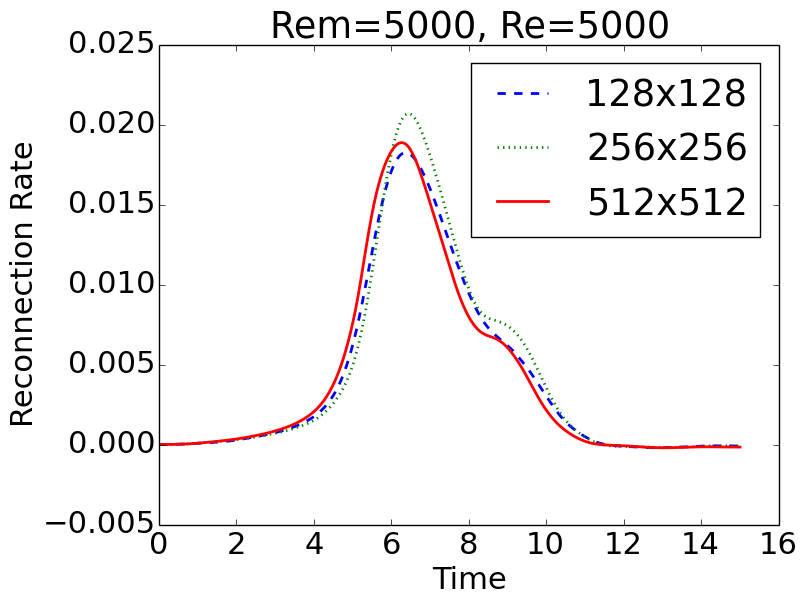} &
 		\includegraphics[width=3.9cm]{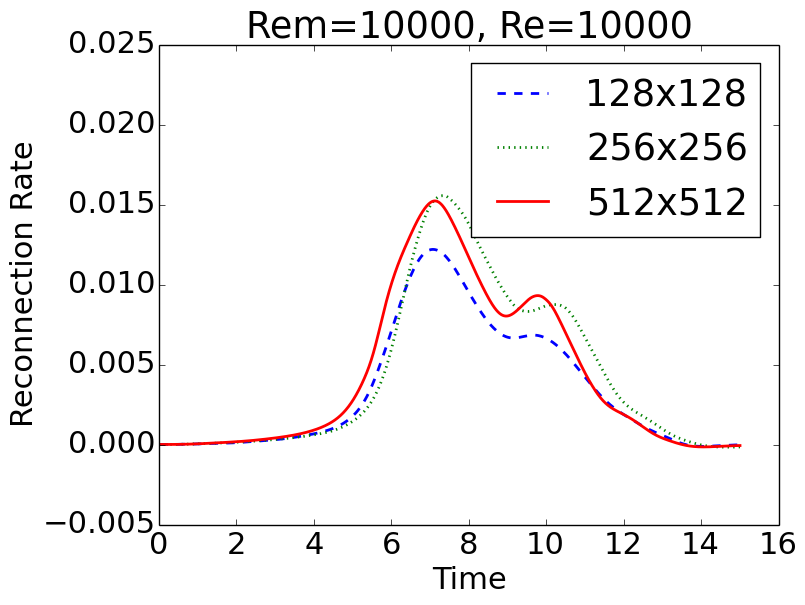} \\
 	\end{tabular}
 	\caption{Reconnection rates for the island coalescence problem.\label{fig:reconnectionrate}}
 \end{figure}
 
 Furthermore, we performed a weak parallel scaling test for nine different combinations of the Reynolds numbers. We chose a coarse grid of $16\times16$ cells with three (1.1M DoFs), four (4.6M DoFs) and five (18.4M DoFs) levels of refinement. All tests were performed with 16 cores per node on 1, 4 and 16 nodes resulting in 16, 64 and 256 cores for the different refinements. We observed (not reported here) that scaling over the nodes with a fixed number of cores per node provides better results than increasing the number of cores per node. This seems to indicate that our code is mainly limited by the memory bandwidth. Furthermore, we ensured that the numbers of cores used in our simulations evenly divide the number of cells in the 16x16 coarse grid to minimize load imbalances. For an optimal scaling of the patch smoother in the multigrid relaxation, the number of patches (and hence vertices) per processor should also be evenly balanced, but this was not implemented.
 
 In Table \ref{tab:islandcoalscaling}, we report the average runtimes per linear iteration rather than the total runtime to take into account that the numbers of linear and nonlinear iterations change slightly between the different refinements. The runtimes only show a slight increase the more cores are used and hence indicate good weak scaling of our method.
 
 As for the lid-driven cavity problem, we observe excellent robustness of the linear and nonlinear iteration counts with respect to the Reynolds numbers. Both linear and nonlinear solvers converge in either 1 or 2 iterations in the investigated ranges of $\Re$ and $\Rem$. We therefore do not report a table here that shows each iteration count.
 

\begin{table}[htbp!]
	\centering
	\resizebox{\textwidth}{!}{%
		\begin{tabular}{r|ccc|ccc|ccc}
			\toprule
			& \multicolumn{3}{c|}{$128\times 128$ on 16 cores} & \multicolumn{3}{c|}{$256\times 256$ on 64 cores} & \multicolumn{3}{c}{$512\times512$ on 256 cores} \\
			\midrule
			$\Rem \backslash\Re$  &1 &     1,000 &    10,000 &1 &     1,000 &    10,000 &1 &     1,000 &    10,000 \\
			\midrule
			1 & 0.13 & 0.12 & 0.13 & 0.14 & 0.14 & 0.14 & 0.17 & 0.17 & 0.17\\
			1,000 & 0.13 & 0.12 & 0.12 & 0.14 & 0.13 & 0.13 & 0.15 & 0.14 & 0.15\\
			10,000 &  0.12 & 0.11 & 0.12 & 0.14 & 0.13 & 0.13 & 0.15 & 0.15 & 0.15\\
			\bottomrule
	\end{tabular}}
	\caption{Average time per linear iteration in minutes for the two-dimensional island coalescence problem. \label{tab:islandcoalscaling}}
\end{table}

\subsection{Three-dimensional results}

In three dimensions, we observe in general that the stationary problems are harder to solve for high parameters than in two dimensions. We believe that the following three points are mainly responsible for this behavior. First of all, the discretization of the electric field changes from a scalar-valued $\mathbb{CG}$-function to a vector-valued $\mathbb{NED}1$-function with tangential boundary conditions. Moreover, the kernel of the term $\vcurl(\u^n \times \B)$ is larger in three dimensions which degrades the performance of the monolithic solver for the $(\E_h, \B_h)$ block for high $\Rem$. Furthermore, the grids we consider are much coarser than in two dimensions because of computational costs.

\subsubsection{Stationary lid-driven cavity problem in three dimensions}\label{sec:stationaryliddrivencavityproblemin3d}

We adapt the two-dimensional lid-driven cavity problem to three dimensions by considering the domain $\Omega=(-1/2, 1/2)^3$ and the boundary conditions $\u=(1,0,0)^\top$ on the boundary $y=0.5$ and $\u=(0,0,0)^\top$ on the other faces. The background magnetic field $\B_0 = (0,1,0)^\top$ determines the boundary conditions for $\B$. We use a coarse grid of $8\times8\times8$ cells with 3 levels of refinement which results in a $64\times64\times64$ grid with 59.1 million DoFs. For the three-dimensional problem, we only investigate the Newton linearization as we have seen in two dimensions that the Newton iteration outperforms the Picard iteration in nearly all cases. The results on the left in Table \ref{tab:ldc3dstationarySRE} show a good control over the linear iteration numbers for the lid-driven cavity problem, where the case of $S=1$ and $\Re=10,000$ seems to be the most challenging case.

On the right of Table \ref{tab:ldc3dstationarySRE} we report a comparison to taking the outer Schur complement that eliminates the $(\E_h, \B_h)$ block.  As mentioned in Section \ref{sec:derivationofblockpreconditioner}, we can clearly see that this choice performs worse for high values of $\Re$ and $S$ where no convergence in 50 linear iterations was reached. We observed similar behavior for unreported experiments on transient and two-dimensional problems.

We do not include a full table for high $\Rem$, as in this case the monolithic multigrid solver cannot deal with the term $\vcurl(\u^n \times \B^n)$ that occurs in the Newton linearization. As in two dimensions, this term is crucial for the convergence of the nonlinear iteration. For Newton, the iteration counts increase very slightly from $\Re=\Rem=1$ by 8.0 Krylov iterations per nonlinear step to 10.0 iterations for $\Rem=500$ and $\Re=1$ and fails to converge for higher $\Rem$. We want to emphasize that in this case the failure of convergence is indeed caused by the inner multigrid method and not by an inaccurate outer Schur complement approximation. To the best of our knowledge, preconditioning methods that robustly treat the vector Laplace operator with an additional $\vcurl(\u^n \times \B)$ term in three dimensions are not known, and we intend to investigate this problem further in future work.

\begin{table}[htbp!]
	\centering
	\begin{minipage}[c]{0.45\textwidth}
		\resizebox{\textwidth}{!}{
			\begin{tabular}{r|ccc}
				\toprule 
				\multicolumn{4}{c}{Using $\tilde{\mathcal{S}}^{(\u, p)}$ for order $(\u, p, \E, \B)$} \\
				\midrule
				$S \backslash\Re$  &1 &     1,000 &    10,000 \\
				\midrule
				1  & ( 3) 6.0 & ( 3) 7.0 & ( 4)20.0 \\
				1,000 & ( 3) 7.3 & ( 2) 9.5 & ( 2) 6.5\\
				10,000 & ( 3) 9.0 & ( 2)13.0   & ( 2)12.5  \\
				\bottomrule
		\end{tabular}}
	\end{minipage}
	\begin{minipage}[c]{0.43\textwidth}
		\resizebox{\textwidth}{!}{
			\begin{tabular}{r|ccc}
				\toprule
					\multicolumn{4}{c}{Using $\tilde{\mathcal{S}}^{(\E, \B)}$ for order $(\E, \B, \u, p)$} \\
				\midrule
				$S \backslash\Re$  &1 &     1,000 &    10,000 \\
				\midrule
				1  &( 3) 6.0 & (4)14.7 & (-)$>$50\\
				1,000 &  ( 3)12.7 & (-)$>$50 &(-)$>$50 \\
				10,000 &  ( 3)20.0 & (-)$>$50 &(-)$>$50\\
				\bottomrule
		\end{tabular}}
	\end{minipage}
	\caption{(left) Iteration counts for the stationary lid-driven cavity problem in 3D for the Newton linearization. (right) Iteration counts for the stationary lid-driven cavity problem in 3D for taking the outer Schur complement that eliminates the $(\E_h, \B_h)$ block.}
	\label{tab:ldc3dstationarySRE}
\end{table}

\subsubsection{Time-dependent lid-driven cavity problem in three dimensions}
Finally, we consider the time-dependent version of the three-dimensional lid-driven cavity problem, which was also investigated in detail in \cite{Ma2016}.
The numerical results in Table \ref{tab:ldc3dtimedep} show good control of the iteration counts and the linear iteration numbers only notably increases for very high values of $S$. Moreover, we observe robust convergence of the monolithic multigrid solver for the $(\E_h, \B_h)$ block for high $\Rem$. As mentioned before in Section \ref{sec:solverformagneticblock}, this can be explained by the fact that the problem does not become nearly singular for high $\Rem$ due to the extra mass matrix. Therefore, the fact that the kernel of $\vcurl(\u^n\times \B)$ is not fully captured by the multigrid method has less influence. 
\newline 

\begin{minipage}[c]{0.45\textwidth}
	  \resizebox{\textwidth}{!}{
    \begin{tabular}{r|ccc}
      \toprule
      $S \backslash\Re$  &1 &     10,000 &    100,000 \\
      \midrule
      1  &   (2.1) 7.3 & (3.2) 2.1 & (3.3) 2.0\\
      1,000 & (3.0) 8.6 & (3.3) 2.8 & (3.5) 2.6 \\
      10,000 & (4.0)11.3 & (4.0) 7.0 & (4.0) 6.2 \\
      \bottomrule
    \end{tabular}}
\end{minipage}
\begin{minipage}[c]{0.463\textwidth}
	 \resizebox{\textwidth}{!}{
      \begin{tabular}{r|ccc}
	\toprule
	$\Rem \backslash\Re$  &1 &     10,000 &    100,000 \\
	\midrule
	1  &     (2.1) 7.3 & (3.2) 2.1 & (3.3) 2.0\\
	10,000 & (2.5) 7.1   & (3.2) 2.0 & (3.3) 2.0   \\
	100,000 & (3.0)15.1 & (3.2) 2.6 & (3.3) 2.0 \\
	\bottomrule
\end{tabular}}
\end{minipage}
\begin{table}[htbp!]
	\centering
	\resizebox{0.45\textwidth}{!}{
		\begin{tabular}{r|ccc}
			\toprule
			$\Rem \backslash S$  &1 &     100 &    1,000 \\
			\midrule
			1  & (2.1) 7.3 & (2.1) 7.3 & (3.0) 8.6\\
			1,000 & (3.0) 7.6 & (3.0) 6.8& (3.1) 9.5\\
			10,000 & (2.5) 7.1  & (3.1) 7.1&(3.2) 9.7\\
			\bottomrule
	\end{tabular}}
\caption{Iteration counts for the transient lid-driven cavity problem in 3D for the Newton linearization.\label{tab:ldc3dtimedep}}
\end{table}
 

\section{Conclusion and outlook}\label{sec:conclusionandoutlook}
We have presented scalable block preconditioners for an augmented Lagrangian formulation of the incompressible MHD equations that exhibit parameter-robust iteration counts in most cases. We described how to control the outer Schur complement of two linearization types and introduced a special monolithic multigrid method to solve the electromagnetic block. This method is fully $\Rem$-robust in two dimensions; in three dimensions, it is able to efficiently compute results for higher parameters than was previously possible. Furthermore, our solvers allow the use of fully implicit methods for time-dependent problems.
We aim to include stabilization techniques for high magnetic Reynolds numbers in future work and further investigate how to develop a robust multigrid method for the problem including the term $\vcurl(\u^n \times \B)$. This would enable a more robust solver for the most difficult case of stationary problems in three dimensions at high magnetic Reynolds numbers.

\section*{Acknowledgement}
The authors would like to thank Kaibo Hu for many useful suggestions and discussions.

\section*{Code availability}
The code that was used to generate the numerical results and all major Firedrake components have been archived on \cite{zenodo/Firedrake-20220119.0}.

\bibliographystyle{siamplain}
\bibliography{literature}

\end{document}